\documentclass[11pt,twoside]{amsart}
\usepackage{amsfonts}
\usepackage{epsfig,graphics}
\usepackage{amssymb}
\usepackage{amscd}

\newtheorem{theoreme}{Th\'eor\`eme}[section]

\newtheorem{proposition}[theoreme]{Proposition}

\newtheorem{definition}[theoreme]{D\'efinition}

\newtheorem{lemme}[theoreme]{Lemme}

\newtheorem{fait}[theoreme]{Fait}

\newtheorem{remarque}[theoreme]{Remarque}

\newtheorem{notation}[theoreme]{Notations}

\newcommand{\R}{{\bf R}}

\newcommand{\N}{{\bf N}}

\newcommand{\RP}{{\bf RP}}

\newcommand{\Ker}{\text{Ker }}

\newcommand{\Hom}{\text{Hom }}
\newcommand{\End}{\text{End }}

\newcommand{\diag}{\text{diag\;}}

\newcommand{\Ein}{\text{Ein}}

\newcommand{\Conf}{{\mathcal Conf }}
\newcommand{\Cont}{{\mathcal C}^0}

\newcommand{\Ad}{\text{Ad }}

\newcommand{\lien}{{\mathfrak{n}}}
\newcommand{\liez}{{\mathfrak{z}}}
\newcommand{\lieg}{{\mathfrak{g}}}
\newcommand{\lieh}{{\mathfrak{h}}}
\newcommand{\liep}{{\mathfrak{p}}}
\newcommand{\lies}{{\mathfrak{s}}}
\newcommand{\lier}{{\mathfrak{r}}}
\newcommand{\oo}{{\mathfrak{o}}}

\newcommand{\om}{\ensuremath{{ \omega^M}}}
\newcommand{\on}{\ensuremath{{ \omega^N}}}

\newcommand{\hm}{\ensuremath{{ \hat{M}}}}
\newcommand{\hn}{\ensuremath{{ \hat{N}}}}
\newcommand{\hx}{\ensuremath{{ \hat{x}}}}
\newcommand{\hy}{\ensuremath{{ \hat{y}}}}
\newcommand{\hz}{\ensuremath{{ \hat{z}}}}

\newenvironment{preuve}{\medskip \noindent {\bf Preuve : }}
   {$\diamondsuit$ }

\begin{document}
\pagenumbering{arabic}
\title[D\'eg\'enerescence  locale]{D\'eg\'enerescence locale des transformations conformes pseudo-riemanniennes}
\author{Charles Frances}
\date{\today}
\address{Charles Frances. 
Laboratoire de
Math\'ematiques, 
Universit\'e Paris-Sud. 
91405 ORSAY Cedex.}
\email{Charles.Frances@math.u-psud.fr}
\maketitle

\begin{abstract}
Nous \'etudions l'ensemble $\Conf (M,N)$ des immersions conformes entre deux vari\'et\'es pseudo-riemanniennes $(M,g)$ et $(N,h)$.  Nous caract\'erisons notamment l'adh\'erence de $\Conf(M,N)$ dans l'espace des fonctions continues ${\mathcal C}^0(M,N)$, et d\'ecrivons quelques propri\'et\'es g\'eom\'e\-triques de $(M,g)$ lorsque cette adh\'erence est non vide.   \end{abstract}

\section{Introduction}

Cet article porte sur l'\'etude de l'ensemble $\Conf(M,N)$ des immersions conformes lisses, entre deux vari\'et\'es connexes pseudo-riemanniennes $(M,g)$ et $(N,h)$, \'egalement lisses, de m\^eme signature $(p,q)$.  Il s'agit des immersions $f:M \to N$ qui satisfont $f^*h=e^{\sigma}g$, avec $\sigma$ une fonction lisse sur $M$. On sera particuli\`erement int\'eress\'e par  les mani\`eres dont peuvent d\'eg\'en\'erer de telles immersions conformes, autrement dit, par ``l'infini" de $\Conf(M,N)$.  Un moyen de comprendre cet infini consiste \`a  d\'ecrire l'adh\'erence de $\Conf(M,N)$ dans ${\mathcal C}^0(M,N)$, l'espace des applications continues de $M$ dans $N$, que l'on munit de la topologie de la convergence uniforme sur les compacts de $M$.

Lorsque l'on est en pr\'esence de deux vari\'et\'es riemanniennes $(M,g)$ et $(N,h)$ de m\^eme dimension $n \geq 2$, on dispose de r\'esultats tr\`es g\'en\'eraux,  prouv\'es par   J. Ferrand dans \cite{ferrand2}, qui traitent de l'ensemble $Q_K(M,N)$ des {\it plongements} $K$-quasi-conformes entre $(M,g)$ et $(N,h)$, o\`u $K \geq 1$  (voir \'egalement les travaux \cite{gehring} et \cite{vassala} dans le cadre de l'espace euclidien).  Il ressort essentiellement des th\'eor\`emes $\mbox{A}$, $\mbox{B}$ et $\mbox{C}$ de \cite{ferrand2},  que si une suite $(f_k)$ de $Q_K(M,N)$ admet une limite $f$ dans ${\mathcal C}^0(M,N)$, alors $f$ appartient \`a $ Q_K(M,N)$ ou $f$ est une application constante.  De plus, si l'on est dans ce dernier cas, la vari\'et\'e $(M,g)$ est $K$-quasi-conforme \`a un ouvert de $\R^n$. Autrement dit, les applications quasi-conformes entres vari\'et\'es riemanniennes d\'eg\'en\`erent comme le font les transformations de M\"obius.

Notre objectif est de comprendre ce qu'il advient en signature quelconque,  pour les {\it immersions} conformes, donc sans hypoth\`ese d'injectivit\'e, mais en se limitant aux applications lisses, et \`a la dimension $\geq 3$.  On travaille donc dans le contexte des structures g\'eom\'etriques rigides.

Notre premier r\'esultat caract\'erise les {\it familles normales} de $\Conf(M,N)$, c'est-\`a-dire celles qui sont relativement compactes dans ${\mathcal C}^0(M,N)$.  Nous allons en fait donner un crit\`ere local simple, qui assure la relative compacit\'e au sein des applications lisses.

\begin{theoreme}
\label{thm.normal}
Soient $(M,[g])$ et $(N,[h])$ deux structures conformes pseudo-riemanniennes de signature $(p,q)$, $p+q \geq 3$.  Soit  ${\mathcal F}$ une famille d'immersions conformes lisses de $(M,[g])$ dans $(N,[h])$.  On suppose qu'il existe $x \in M$ tel que~:
\begin{enumerate}
\item{L'ensemble $E:=\{f(x) \ | \ f \in {\mathcal F}   \}$ est relativement compact dans $N$.}
\item{La famille ${\mathcal F}$ est \'equicontinue en $x$.}
\end{enumerate}
Alors il existe un ouvert $U \subset M$ contenant $x$ tel que ${\mathcal F}_{|U}$ soit relativement compacte dans ${\mathcal C}^{\infty}(U,N)$, pour la topologie de la convergence $C^{\infty}$ sur les compacts de $U$.
\end{theoreme}

%\begin{corollaire}
%\label{coro.regularite}
%Si une suite d'immersions conformes $f_k : (M,g) \to (N,h)$ 
%\end{corollaire}

On veut maintenant d\'ecrire quelles  applications  peuvent appara\^{\i}tre dans l'adh\'erence de $\Conf(M,N)$ dans ${\mathcal C}^0(M,N)$.  Le th\'eor\`eme \ref{thm.normal} assure d\'ej\`a que ces applications sont n\'ecessairement lisses.  On peut en fait dire beaucoup plus gr\^ace au th\'eor\`eme ci-dessous. Nous rappelons qu'une sous-vari\'et\'e $\Sigma$ de $(M,g)$ est dite d\'eg\'en\'er\'ee lorsque la restriction de $g$ \`a $\Sigma$ est d\'eg\'en\'er\'ee. Si $x \in \Sigma$, le radical isotrope en $x$ est le plus grand sous-espace de $T_x\Sigma$ sur lequel la forme $g_x$ s'annule. La vari\'et\'e $\Sigma$ est dite totalement isotrope lorsque la restriction de $g$ \`a $\Sigma$ est identiquement nulle. Dans le cadre riemannien, les sous-vari\'et\'es d\'eg\'en\'er\'ees sont les points. Enfin, la notion de sous-vari\'et\'e totalement g\'eod\'esique conforme sera introduite en Section \ref{sec.geodesiques}.

\begin{theoreme}
\label{thm.applications.limites}
Soient $ (M,[g])$ et $ (N,[h])$ deux structures conformes  pseudo-riemanniennes de signature $(p,q)$ , $p+q\geq 3$. Soit $ f_k : (M,[g]) \to (N,[h])$ une suite d'immersions conformes lisses qui converge uniform\'ement sur les compacts de $M$  vers une application $f \in {\mathcal C}^0(M,N)$.  Alors l'application $f$ est de classe ${\mathcal C}^{\infty}$ et de rang constant, et la convergence de $(f_k)$ vers $f$ est ${\mathcal C}^{\infty}$ sur les compacts de $M$. De plus, on est dans exactement l'un des trois cas suivants~:
\begin{enumerate}
\item{L'application $f$ appartient \`a $\Conf (M,N)$.}
\item{L'application $f$ est constante.  Dans ce cas, $(M,g)$ est localement conform\'ement Ricci-plate.}
\item{
\begin{enumerate}
\item{Localement, $f$ est une submersion sur une sous-vari\'et\'e lisse, totalement isotrope de $N$.} 
\item{Les fibres de $f$ sont des sous-vari\'et\'es totalement g\'eod\'esiques conformes, qui sont d\'eg\'en\'er\'ees. Le radical isotrope de ces fibres a pour dimension le rang de $f$.}

\end{enumerate}}
\end{enumerate}
\end{theoreme}

Dans le cadre riemannien, le cas $(3)$ ne peut pas se produire. On retrouve que la limite d'une suite d'immersions conformes riemanniennes est soit une immersion conforme, soit une application constante.   En signature lorentzienne, les limites possibles de suites d'immersions conformes sont soit des immersions conformes, soit des applications constantes, soit des submersions sur une g\'eod\'esique de lumi\`ere de $(N,h)$, dont les fibres sont des hypersurfaces d\'eg\'en\'er\'ees, totalement g\'eod\'esiques conformes, de $(M,g)$. 

Dans tous les cas,  l'\'enonc\'e pr\'ec\'edent montre  que lorsque $\Conf(M,N)$ n'est pas ferm\'e dans ${\mathcal C}^0(M,N)$, la vari\'et\'e $(M,g)$ pr\'esente des propri\'et\'es g\'eom\'e\-triques int\'eressantes. Dans les cadres riemanniens et lorentziens, on peut pr\'eciser ces r\'esultats  pour obtenir~:
\begin{theoreme}
\label{thm.ferme.ou.plat}
Soient $ (M,[g])$ et $ (N,[h])$ deux structures pseudo-rieman\-niennes de m\^eme  signature $(p,q)$, $p+q \geq 3$.  On suppose que $\Conf (M,N)$ n'est pas ferm\'e dans ${\mathcal C}^0(M,N)$. Alors~:
\begin{enumerate}
\item{ Si $(M,[g])$ est riemannienne, elle est localement  conform\'ement plate.}
\item{Si $(M,[g])$ est lorentzienne, alors elle est localement conform\'ement Ricci-plate.  Si de plus l'adh\'erence de $\Conf (M,N)$  dans ${\mathcal C}^0(M,N)$ contient une application constante, alors $(M,[g])$ est localement conform\'ement plate.}
\end{enumerate}
\end{theoreme}

Enfin, notre dernier th\'eor\`eme est de nature dynamique.  \'Etant donn\'ee une suite $(f_k)$ de $\Conf (M,N)$ qui converge dans ${\mathcal C}^0(M,N)$ vers une application $f$, le Th\'eor\`eme \ref{thm.applications.limites}  nous dit que les fibres de $f$ sont des sous-vari\'et\'es de $M$.  Si $x$ et $y$ sont dans la m\^eme fibre, alors $f_k(x)$ et $f_k(y)$ ont bien entendu la m\^eme limite.  La question est maintenant de savoir ``\`a quelle vitesse" les points $f_k(x)$ et $f_k(y)$ se rapprochent. Le r\'esultat qui suit peut \^etre vu comme un pendant  conforme de \cite[Th\'eor\`emes 1.1 et 1.2]{ghani1}, et donne l'existence de  ``vitesses critiques" naturellement associ\'ees \`a la suite $(f_k)$, qui d\'efinissent \`a leur tour des sous-feuilletages  totalement g\'eod\'esiques conformes des fibres de $f$.  Avant de donner l'\'enonc\'e pr\'ecis, rappelons que si $(a_k)$ et $(b_k)$ sont deux suites de r\'eels positifs, alors $a_k=\Theta(b_k)$ signifie qu'il existe deux constantes $C_1,C_2>0$ telles que pour $k$ suffisamment grand~:
$$ C_1b_k \leq a_k \leq C_2 b_k.$$
Par $a_k=O(b_k)$, on entend qu'il existe $C>0$ tel que pour $k$ suffisamment grand, $a_k \leq C b_k$.

\begin{theoreme}
\label{thm.integrabilite-stable}
Soit $(f_k)$ une suite de $\Conf(M,N)$ qui converge uniform\'e\-ment sur les compacts de $M$ vers une application  limite  $f$.  Quitte \`a remplacer $(f_k)$ par une suite extraite, il existe~:
\begin{itemize}
\item{ Un entier $s \geq 1$ et une filtration de $TM$, ${\mathcal F}_{0}=\{ 0 \} \subsetneq {\mathcal F}_{1} \subsetneq \ldots  \subsetneq {\mathcal F}_{s-1} \subsetneq TM$, qui s'int\`egre en $s$ feuilletages lisses totalement g\'eod\'esiques conformes de $M$, $F_{0} \subsetneq \ldots \subsetneq F_{s-1}$.}
\item{ Des suites convergentes $\mu_1(k), \ldots, \mu_s(k)$ de $\R_+^*$, v\'erifiant  $\mu_j(k)=\mbox{O}(\mu_{j+1}(k))$ pour $j \in \{1, \ldots, s-1  \}$, et des entiers $n_1, \ldots , n_s$ de $\N^*$ avec $n_1+ \ldots +n_s=n$, tels que pour tout $k \in \N$, la matrice  $$   \left( \begin{array}{ccc}
\mu_1(k)I_{n_1} &  & 0 \\
  & \ddots   &  \\
  0&    & \mu_s(k)I_{n_s}
\end{array} \right)$$ pr\'eserve la classe conforme de la forme quadratique $2x_1x_{p+q}+\ldots+2x_px_{q+1}+\Sigma_{p+1}^{q}x_i^2$,}
 \end{itemize}
satisfaisant les propri\'et\'es suivantes~:
\begin{enumerate}
%\item{Pour chaque $j \in \{-s, \ldots, s \}$ la suite $\mu_j(k)$ a une limite dans $\R$.}
%\item{Pour tout $j \in \{ -s, \ldots , s-1 \}$, $\frac{\mu_{j+1}(k)}{\mu_{j}(k)} \to \infty$.}
%\item{Pour tout $x \in U$, ${\mathcal F}_i(x) \bot {\mathcal F}_j(x)$ si $i \not = -j$.  Par ailleurs la somme restriction de $g$ \`a ${\mathcal F}_{-i}(x) \oplus {\mathcal F}_i(x)$ n'est pas  d\'eg\'en\'er\'ee.}
\item{ 
\begin{enumerate}
\item{Un vecteur non nul $u \in T_xM$ appartient \`a  $T_xM \setminus {\mathcal F}_{s-1}(x)$, si et seulement si   pour toute suite $(u_k)$ de $T_xM$ qui converge vers $u$, 
$$||D_xf_k(u_k)|| = \Theta(  \mu_s(k)).$$}
\item{   Un vecteur non nul $u \in T_xM$ appartient \`a ${\mathcal F}_{j}(x) \setminus {\mathcal F}_{j-1}(x)$, $j=1, \ldots s-1$, si et seulement si les deux conditions ci-dessous sont satisfaites~:
\begin{enumerate}
\item{Pour toute suite $(u_k)$ de $T_xM$ qui converge vers $u$,  
$$ \mu_{j}(k)=O(||D_xf_k(u_k)|).$$} 
\item{   Il existe une suite $(u_k)$ de $T_xM$ qui converge vers $u$ telle  que $ ||D_xf_k(u_k)||  = \Theta(\mu_j(k))$.}
\end{enumerate}}
\end{enumerate} }
\item{Chaque $x \in M$ admet un voisinage $U_x$ tel que la feuille locale $F_j^{loc}(x)$ de $x$ dans $U_x$ soit caract\'eris\'ee par :

\begin{enumerate}
\item{Un point $y$  appartient \`a $U_x \setminus  F_{s-1}^{loc}(x)$ si et seulement si pour toute suite $(y_k)$ de $U_x$ qui converge vers $y$,  $d(f_k(x),f_k(y_k)) = \Theta(\mu_s(k))$.}

\item{   Un point $y$ appartient \`a  $ F_{j}^{loc}(x) \setminus F_{j-1}^{loc}(x)$, $j=1, \ldots s-1$, si et seulement si les deux conditions ci-dessous sont satisfaites~:
\begin{enumerate} 
\item{Pour toute suite $(y_k)$ de $U_x$ qui converge vers $y$, 
$$  \mu_{j}(k)=O(d(f_k(x),f_k(y_k))).$$} 
 \item{Il existe une suite $(y_k)$ de $U_x$ qui converge vers $y$ telle que $d(f_k(x),f_k(y_k)) = \Theta(\mu_j(k))$.}
 \end{enumerate}}
 \end{enumerate} }
\end{enumerate}

\end{theoreme}

Nous verrons en Section \ref{sec.unicite} que l'entier $s$, les sous-vari\'et\'es ${ F}_j^{loc}(x)$ et les suites $\mu_j(k)$, s'ils v\'erifient les conclusions $2 \ (a),(b)$ du th\'eor\`eme, sont uniques.

\section{Pr\'erequis alg\'ebriques et g\'eom\'etriques}
\label{sec.geometrie}
%\subsection{G\'eom\'etries de Cartan}
 
\subsection{L'univers d'Einstein}
\label{sec.univers}
%Dans de nombreux contextes g\'eom\'etriques, il existe  un objet  privil\'egi\'e, souvent caract\'eris\'e par le fait qu'il admet beaucoup de sym\'etries.  Dans le cadre  des structures conformes de signature $(p,q)$, cet objet est l'univers d'Einstein $\Ein^{p,q}$, dont nous donnons maintenant une br\`eve description  (voir \cite{primer} ou \cite{charlesthese} pour une \'etude plus d\'etaill\'ee).  Nous faisons dans tout l'article la convention $p \leq q$ (dans la cas d'une signature o\`u $q>p$, il faut permuter les mots ``temps" et ``espace" dans tous les \'enonc\'es de l'article).
Nous commen\c{c}ons par de brefs rappels g\'eom\'e\-triques sur l'espace mod\`ele conform\'ement plat de signature $(p,q)$. Pour une \'etude plus d\'etaill\'ee, nous renvoyons le lecteur \`a  \cite{primer} ou \cite{charlesthese}.

Soit $\R^{p+1,q+1}$  l'espace $\R^{p+q+2}$ muni de la forme quadratique :
$$ Q^{p+1,q+1}(x):=2x_0x_{p+q+1}+...+2x_px_{q+1}+\Sigma_{p+1}^q x_i^2$$
Le c\^one isotrope de cette forme quadratique est not\'e ${\mathcal N}^{p+1,q+1}$.  La restriction de $Q^{p+1,q+1}$ fournit une m\'etrique d\'eg\'en\'er\'ee sur ${\mathcal N}^{p+1,q+1} \setminus \{0 \}$, dont le noyau est de dimension $1$, et est tangent aux g\'en\'eratrices du c\^one ${\mathcal N}^{p+1,q+1} $.  Ainsi, le projectivis\'e ${\bf P}({\mathcal N}^{p+1,q+1} \setminus \{ 0 \})$ est une sous-vari\'et\'e lisse de $\RP^{p+q+1}$, naturellement munie d'une classe conforme de m\'etriques non d\'eg\'en\'er\'ees, de signature $(p,q)$.  
On appelle {\it univers d'Einstein} de signature $(p,q)$, et l'on note $\Ein^{p,q}$, cette vari\'et\'e compacte ${\bf P}({\mathcal N}^{p+1,q+1} \setminus \{  0  \})$ munie de la structure conforme ci-dessus.

L'espace $\Ein^{0,q}$ n'est autre que la sph\`ere ${\bf S}^q$ munie de la classe conforme de la m\'etrique ``ronde" $g_{{\bf S}^q}$.
Pour $p \geq 1$,  le produit ${\bf S}^p \times {\bf S}^q$, muni de la classe conforme de la m\'etrique produit $-g_{{ \bf S}^p} \oplus g_{{\bf S}^q}$, est  un rev\^etement double, conforme, de $\Ein^{p,q}$.

Soit $\mbox{O}(p+1,q+1)$ le groupe des transformations lin\'eaires qui laissent $Q^{p+1,q+1}$ invariante.  L'action naturelle de $\mbox{O}(p+1,q+1)$ sur $\Ein^{p,q}$ pr\'eserve la structure conforme de $\Ein^{p,q}$, et il s'av\`ere que $\mbox{O}(p+1,q+1)$ est tout le groupe des diff\'eomorphismes conformes de $\Ein^{p,q}$  (voir \cite{charlesthese} ou \cite{schot}, ainsi que les r\'ef\'erences de \cite{kuehnel3}).

Dans tout l'article nous appellerons $o$ le point de $\Ein^{p,q}$ correspondant \`a $[e_0]$.  Son stabilisateur $P \subset \mbox{O}(p+1,q+1)$ est un sous-groupe parabolique isomorphe au produit semi-direct $(\R_+^* \times \mbox{O}(p,q)) \ltimes \R^{p,q}$, avec $n=p+q$.  Du point de vue conforme, $\Ein^{p,q}$ est donc l'espace homog\`ene $\mbox{O}(p+1,q+1)/P$.

\subsection{L'alg\`ebre de Lie $\lieg=\oo(p+1,q+1)$}
\label{sec.algebre}
L'alg\`ebre de Lie $\oo(p+1,q+1)$ est compos\'ee des matrices $X$ de taille
$(p+q+2) \times (p+q+2)$ qui satisfont l'identit\'e : 
$$ X^t J_{p+1,q+1} + J_{p+1,q+1} X = 0.$$
Ici,  $J_{p+1,q+1}$ est la matrice dans la base $(e_0, \ldots ,e_{n+1})$ de la forme quadratique $Q^{p+1,q+1}$.

L'alg\`ebre $\oo(p+1,q+1)$ s'\'ecrit comme une somme $\lien^- \oplus \liez \oplus \lies \oplus \lien^+$, o\`u :

$$ \liez =  \left\{ \left( \begin{array}{ccc}
a &  & 0 \\
  & 0   &  \\
  &    & -a
\end{array} \right) \ :    
\qquad 
\begin{array}{c}
 a \in \R  
  
    \end{array} 
\right\}
$$

$$ \lies =  \left\{ \left( \begin{array}{ccc}
0 &  & 0 \\
  & M   &  \\
  &    & 0
\end{array} \right) \ :    
\qquad 
\begin{array}{c}
 
  M \in \oo(p,q) 
    \end{array} 
\right\}
$$

$$ \lien^+= \left\{ \left( \begin{array}{ccc}
0& -x^t.J_{p,q}   &  0\\
  & 0  & x \\
  &   & 0
\end{array} \right) \ :  
\qquad 
\begin{array}{c}
  x\in \R^{p,q} 
\end{array}
\right\}
$$

$$ \lien^- = \left\{ \left( \begin{array}{ccc}
0&    &  \\
 x & 0  &  \\
 0 & -x^t.J_{p,q}  & 0
\end{array} \right) \ :  
\qquad 
\begin{array}{c}
  x\in \R^{p,q} 
\end{array}
\right\}
$$

On notera \'egalement $\lier=\liez \oplus \lies$.  Par $\mbox{N}^+,$ $\mbox{N}^-$ et $\mbox{Z}$, on d\'esignera les sous-groupes connexes de $\mbox{O}(p+1,q+1)$ dont l'al\`ebre de Lie est $\lien^+$, $\lien^-$ et $\liez$ respectivement.  Il existe par ailleurs dans $\mbox{O}(p+1,q+1)$ un sous-groupe isomorphe \`a $\mbox{O}(p,q)$, dont l'alg\`ebre de Lie est $\lies$~: il sera not\'e  $\mbox{S}$.  Enfin, on appelera $R$ le produit $\mbox{Z} \times \mbox{S}$, dont l'alg\`ebre de Lie est $\lier$.  Ce groupe est isomorphe au produit $\R_+^* \times \mbox{O}(p,q)$.

Dans $\oo(p+1,q+1)$, on appelle ${\mathfrak a}$ l'alg\`ebre constitu\'ee des matrices :
$${\mathfrak a} = \left\{ \left( \begin{array}{ccccccc}
\alpha_1&    & & & & &  \\ $$
 & \ddots & & & & &  \\
 & & \alpha_{p+1} & & &  &\\
 & & & 0_{q-p} & & &  \\
 & & & & -\alpha_{p+1}& & \\
 & & & & & \ddots & \\
 & & & & & & -\alpha_{1} \\
 \end{array}
\right)  :  
\qquad 
\begin{array}{c}
  \alpha_1, \ldots ,\alpha_{p+1} \in \R 
\end{array} \right\}.
$$

Dans l'expression ci-dessus, $0_{q-p}$ d\'esigne la matrice nulle de taille $q-p$.  On appelle ${\mathfrak a}^+$ le sous-ensemble de ${\mathfrak a}$ pour lequel $\alpha_1 \geq \ldots \geq  \alpha_{p+1} \geq 0$, et l'on pose $\mbox{A}^+:= e^{{\mathfrak a}^+}$  (l'image de ${\mathfrak a}^+$ par l'exponentielle dans $\mbox{O}(p+1,q+1)$).  

\subsection{Le groupe parabolique $\mbox{P}$}
\label{sec.matriciel}
Le groupe $\mbox{P}$ est le stabilisateur du point $o=[e_0]$ dans $\mbox{O}(p+1,q+1)$.  Il s'agit du sous-groupe $(\mbox{Z} \times \mbox{S}) \ltimes {\mbox{N}}^+ \subset \mbox{O}(p+1,q+1)$. Comme nous l'avons vu, il est isomorphe au produit semi-direct $(\R_+^* \times \mbox{O}(p,q)) \ltimes \R^{p,q}$, {\it i.e} au groupe conforme de l'espace $\R^{p,q}$. Cet isomorphisme n'est pas fortuit.  Il existe en effet une immersion conforme $j: \R^{p,q} \to \Ein^{p,q}$, appel\'ee {\it projection st\'er\'eographique}, et donn\'ee en coordonn\'ees projectives par la formule $x=(x_1,\ldots,x_n) \mapsto [-\frac{Q^{p,q}(x)}{2},x_1,\ldots,x_2,1]$.  L'image  $j(\R^{p,q})$ est le compl\'ementaire dans $\Ein^{p,q}$ du c\^one de lumi\`ere passant par $o$.  

\subsubsection{Le groupe $\mbox{P}$ comme groupe de transformations affines}
L'application $j$ conjugue l'action naturelle de $\mbox{P}$ sur le compl\'ementaire du c\^one de lumi\`ere issu de $o$, \`a l'action affine de $(\R_+^* \times \mbox{O}(p,q)) \ltimes \R^{p,q}$ sur $\R^{p,q}$.  Gr\^ace \`a la projection st\'er\'eographique $j$, on verra souvent les \'el\'ements de $\mbox{P}$ {\it comme des transformations  affines} de $\R^{p,q}$, not\'ees $A+T$ , o\`u $A \in \R_+^* \times \mbox{O}(p,q)$ est la partie lin\'eaire, et $T \in \R^{p,q}$  le facteur translation.

Vu comme \'el\'ement de $\mbox{O}(p+1,q+1)$, une  translation de vecteur $v \in \R^{p,q}$ admet l'expression~:
\begin{equation}
\label{eq.matricielle}
 n^+(v):=\left( \begin{array}{ccc}
1& -v^t.J_{p,q}   &  -\frac{Q^{p,q}(v)}{2}\\
  & 1  & v \\
  &   & 1
\end{array} \right).  
\end{equation}
L'ensemble des translations forme  le groupe $N^+$, d'alg\`ebre de Lie $\lien^+$.  L'appli\-cation $n^+ : \R^{p,q} \to N^+$ est un isomorphisme de groupes.

Vu dans $\mbox{O}(p+1,q+1)$, un \'el\'ement $\lambda A \in \R_+^* \times \mbox{O}(p,q)$ s'exprime comme~:
\begin{equation}
\label{eq.matricielle2}
 \left( \begin{array}{ccc}
\lambda & 0 & 0 \\
0  & A   & 0 \\
0  & 0   & \frac{1}{\lambda}
\end{array} \right).
\end{equation}

Par ces identifications, un  \'el\'ement $h=\left( \begin{array}{ccccccc}
\alpha_1&    & & & & &  \\ 
 & \ddots & & & & &  \\
 & & \alpha_{p+1} & & &  &\\
 & & & I_{q-p} & & &  \\
 & & & & \frac{1}{\alpha_{p+1}}& & \\
 & & & & & \ddots & \\
 & & & & & & \frac{1}{\alpha_{1}} \\
 \end{array}
\right) $ de $\mbox{A}^+$ est la transformation lin\'eaire diagonale~:
$$h=\alpha_1\left( \begin{array}{ccccccc}
\alpha_2&    & & & & &  \\ 
 & \ddots & & & & &  \\
 & & \alpha_{p+1} & & &  &\\
 & & & I_{q-p} & & &  \\
 & & & & \frac{1}{\alpha_{p+1}}& & \\
 & & & & & \ddots & \\
 & & & & & & \frac{1}{\alpha_{2}} \\
 \end{array}
\right), $$ agissant sur  $\R^{p,q}$.  On note~: $h=\diag(\lambda_1, \ldots, \lambda_n)$, avec $\lambda_1 \geq  \ldots \geq \lambda_n \geq1$.
Il est \`a noter que l'action adjointe  $(\Ad h)$ restreinte \`a $\lien^-$ se fait via la transformation $\diag(\frac{1}{\lambda_1}, \ldots , \frac{1}{\lambda_n})$.

Dans tout l'article, on adoptera tant\^ot la notation $h \in R$, lorsque $h$ est vu comme un \'el\'ement de $\mbox{O}(p+1,q+1)$, tant\^ot $h \in \R_+^* \times \mbox{O}(p,q)$ pour signifier que l'on voit $h$ comme une transformation lin\'eaire de l'espace $\R^{p,q}$.

\section{G\'eom\'etrie sur le fibr\'e de Cartan associ\'e \`a une structure conforme}
Nous commen\c{c}ons par rappeler l'interpr\'etation des structures conformes pseudo-riemanniennes de dimension sup\'erieure ou \'egale \`a 3 en termes de g\'eom\'etrie de Cartan. 

\subsection{Le probl\`eme d'\'equivalence}
\label{sec.principe-equivalence}
Soit $G$ un groupe de Lie, $P \subset G$ un sous-groupe ferm\'e, et ${\bf X}=G/P$.  On appelle {\it g\'eom\'etrie de Cartan } model\'ee sur ${\bf X}$ la donn\'ee d'un triplet $(M,\hm,\omega)$, o\`u :
\begin{enumerate}
\item{$M$ est une vari\'et\'e de m\^eme dimension que ${\bf X}$.}
\item{$\pi_M:\hm \to M$ est un $\mbox{P}$ fibr\'e principal au-dessus de $M$.}
\item{la forme $\omega$ est une $1$-forme \`a valeurs dans $\lieg$, satisfaisant les propri\'et\'es suivantes :
\begin{itemize}
\item{ Pour chaque $\hx \in \hm$, $\omega_{\hx} : T_{\hx}\hm \to \lieg$ est un isomorphisme d'espaces vectoriels. }
\item{Pour tout $X \in \lieg$ et $\hx \in \hm$, $\omega_{\hx}(\frac{d}{dt}_{t =0}R_{e^{tX}}.\hx)=X$.}
\item{Pour tout $p \in \mbox{P}$, $(R_p)^*\omega = (\Ad p^{-1}).\omega$.}
\end{itemize}}
\end{enumerate}
Ici, $R_p$ d\'esigne l'action \`a droite de $p \in \mbox{P}$ sur $\hm$, et $e^{tX}$ est l'exponentielle dans le groupe $G$. Une $1$-forme $\omega$ comme ci-dessus s'appelle une {\it connexion de Cartan} sur $\hm$.

Il faut penser \`a une g\'eom\'etrie de Cartan comme un analogue courbe du {\it mod\`ele  plat} $({\bf X},G,\omega_G)$ o\`u $\omega_G$ est la forme de Maurer-Cartan sur $G$.

Prenons \`a pr\'esent l'exemple de l'espace homog\`ene ${\bf X}:=\Ein^{p,q}=\mbox{O}(p+1,q+1)/P$, et d'une g\'eom\'etrie de Cartan $(M,\hm,\omega)$ model\'ee sur $\Ein^{p,q}$.  Pour tout $x \in M$ et $\hx \in \hm$ au-dessus de $x$, il existe un isomorphisme naturel :
$$ \iota_{\hx} : \lieg/\liep \to T_xM$$ d\'efini par $\iota_{\hx}(\overline \xi):=D_{\hx}\pi_M(\omega_{\hx}^{-1}(\xi))$, o\`u $\xi$ est un repr\'esentant quelconque dans $\lieg$ de la classe $\overline \xi \in \lieg/\liep$.  L'isomorphisme $\iota_{\hx}$ satisfait la relation d'\'equivariance~:
\begin{equation}
\label{equ.equivariance}
 \iota_{\hx.p^{-1}}((\Ad p).\overline \xi)= \iota_{\hx}(\overline \xi), \ \forall p \in \mbox{P}.
\end{equation}

Soit $\lambda^{p,q}$ une m\'etrique de signature $(p,q)$ sur $\lien^-$, qui soit $(\Ad S)$-invariante (notons qu'une telle m\'etrique est unique \`a multiplication par un scalaire pr\`es).  Alors ${\mathcal C}=[\lambda^{p,q}]$ est l'unique classe conforme de produits scalaires de signature $(p,q)$ qui soit $(\Ad P)$-invariante sur $\lieg/\liep$, $\iota_{\hx}({\mathcal C})$ d\'efinit une classe conforme de signature $(p,q)$ sur $T_xM$, ind\'ependante du choix de $\hx \in \hm$ au-dessus de $x$.  Autrement-dit, une g\'eom\'etrie de Cartan $(M,\hm,\omega)$ model\'ee sur $\Ein^{p,q}$ d\'efinit une classe conforme $[g]$ de m\'etriques de signature $(p,q)$ sur $M$.  

Le fibr\'e $\hm$ \'etant fix\'e, il existe {\it a priori} beaucoup de connexions de Cartan $\omega$ d\'efinissant la classe $[g]$ sur $M$.  De m\^eme que dans le contexte m\'etrique, il existe une unique connexion compatible sans torsion (la connexion de Levi-Civita), il existe un choix de normalisation ad\'equat qui rend $\omega$ unique (voir la Section \ref{sec.courbure} ci-dessous).  Pr\'ecis\'ement, voir \cite[chap. 7,  Prop. 3.1 p.285]{sharpe}, \cite{kobayashi}, on peut \'enoncer le~:

\begin{theoreme}[E. Cartan]
\label{thm.principe-equivalence}
Soit $(M,[g])$ une structure conforme pseudo-riemannienne de signature $(p,q)$, $p+q \geq 3$. Alors il existe une unique  g\'eom\'etrie de Cartan $(M,\hm,\omega)$ normale, model\'ee sur $\Ein^{p,q}$, d\'efinissant la structure conforme $(M,[g])$ par la proc\'edure ci-dessus. En particulier, tout diff\'eomorphisme conforme local $\phi$ sur $M$ se remonte en un automorphisme local de fibr\'e (renot\'e $\phi$) qui pr\'eserve $\omega$.
\end{theoreme}

Dans la suite, si $(M,[g])$ est une structure conforme pseudo-riemannienne de dimension $n \geq 3$, nous appellerons le triplet $(M,\hm,\omega)$ donn\'e par le th\'eor\`eme \ref{thm.principe-equivalence}, le {\it fibr\'e normal de Cartan} d\'efini par $(M,[g])$. 

%\section{Holonomie}

\subsubsection{Courbure conforme et connexion normale}
\label{sec.courbure}
Pour une g\'eom\'etrie de Cartan $(M,\hm,\omega)$ model\'ee sur un espace ${\bf X}=G/P$, on a une notion de courbure $K$ d\'efinie, pour tout couple de champs de vecteurs $\hat X$ et $\hat Y$ sur $\hm$ par (voir \cite{sharpe} p.176, p. 191)~:
$$ K(\hat X,\hat Y):=d\omega(\hat X, \hat Y)+[\omega(\hat X),\omega(\hat Y)].$$
 En particulier, si $\hat X$ et $\hat Y$ sont $\omega$-constants, on obtient~:
 \begin{equation}
 \label{eq.courbure}
  K(\hat X,\hat Y)=[\omega(\hat X),\omega(\hat Y)]-\omega([\hat X,\hat Y]).
  \end{equation} 
Par ailleurs, en un point o\`u $\hat X$ ou $\hat Y$ est tangent aux fibres de $\hm$, la courbure s'annule.  
On peut donc \'egalement voir la courbure comme une application $\kappa: \hm \to \Hom(\Lambda^2(\lieg/\liep),\lieg)$. Cette fonction courbure s'annule au-dessus d'un ouvert $U \subset M$ si et seulement si cet ouvert est localement conform\'ement plat. 

L'application $\kappa$ satisfait la relation d'\'equivariance~:
\begin{equation}
\label{equ.equivariance.courbure}
(\Ad p^{-1}).\kappa_{\hx.p^{-1}}((\Ad p).\xi, (\Ad p).\eta)=\kappa_{\hx}(\xi,\eta).
\end{equation}

Dans le cas o\`u $(M,\hm,\omega)$ est le fibr\'e normal de Cartan d'une structure conforme pseudo-riemannienne $(M,[g])$, alors la condition de normalisation sur  $\omega$ dit deux choses~: d'une part que la fonction courbure $\kappa$, au lieu d'\^etre \`a valeurs dans $\Hom(\Lambda^2(\lieg/\liep),\lieg)$ est \`a valeurs dans $\Hom(\Lambda^2(\lieg/\liep),\lies \oplus \lien^+ )$, et d'autre part que la composante $\kappa_{\lies}$ est dans le noyau de l'homomorphisme de Ricci (voir \cite{sharpe} p.280, et la Prop. 3.1 p.285).     Cette composante $\kappa_{\lies}$  de $\kappa$ sur $\lies$ correspond au tenseur de Weyl $W$ sur la vari\'et\'e $M$ (\cite{sharpe} p.236 et p.290), via la formule~:
$W_y(u,v)w=[\kappa_{\lies}(\iota_{\hy}^{-1}(u),\iota_{\hy}^{-1}(v)),\iota_{\hy}^{-1}(w)]$.

\subsection{Application exponentielle conforme}
  Le choix d'un  \'el\'ement $Z$ dans $\lieg=\oo(p+1,q+1)$ d\'efinit naturellement un champ de vecteurs $\hat Z$ sur $\hm$ par la relation $\omega(\hat Z)=Z$.  Si $Z \in \lieg$, on note $\psi_Z^t$ le flot local engendr\'e sur $\hm$ par le champ $\hat Z$.  En chaque $\hx \in \hm$, on d\'efinit   ${\mathcal W}_{\hx} \subset \lieg$ l'ensemble des $Z$ tels que $\psi_Z^t$ est d\'efini pour  $t \in [0,1]$ en $\hx$.  On d\'efinit  l'application exponentielle en $\hx$~:
  $$\exp(\hx, \ ) : {\mathcal W}_{\hx} \to \hm$$
    par :
  $$ \exp(\hx,Z):=\psi_{\hat Z}^1.\hx$$
  
  Il est standard de montrer que ${\mathcal W}_{\hx}$ est un voisinage de $0$, et que l'application $\xi \mapsto \exp(\hx, \xi)$ r\'ealise un diff\'eomorphisme d'un voisinage ${\mathcal V}_{\hx} \subset {\mathcal W}_{\hx}$ contenant $0$ sur un voisinage de $\hx$ dans $\hm$.
 
% On v\'erifie ais\'ement que si $x \in M$, $\hx \in \hm$ au-dessu de $M$, $\xi \in \lieg$ et $\alpha(s):=\pi(\exp(b,s\xi))$, alors le d\'eveloppement $\beta(s):={\mathcal D}_{x}^{b}(\alpha)(s)$ est donn\'e par $\beta(s)=\pi_G(e^{s \xi})$.
 \subsubsection{Application exponentielle et immersions conformes} 
 Soient $(M,[g])$ et $(N,[h])$ deux structures pseudo-riemanniennes  de signature $(p,q)$, $p+q \geq 3$. On appelle $(M,\hm,\om)$ et $(N,\hn,\on)$ les fibr\'es normaux de Cartan associ\'es.  Pour ne pas alourdir les notations, on appellera indistinctement $\exp$ les applications exponentielles sur $\hm$ et  $\hn$.
  Soit $f$ est une immersion conforme de $(M,[g])$ dans $(N,[h])$.  Alors $f$ se remonte  en une immersion, renot\'ee $f$,  du fibr\'e $(M,\hm,\omega)$ dans le fibr\'e $(N,\hn,\on)$, qui satisfait de plus $f^*(\on)=\om$. Soit $Z \in \lieg$, et $\hat Z_1$  (resp. $\hat Z_2$) le champ de vecteurs sur $\hm$  (resp. sur $\hn$) satisfaisant $\om(\hat Z_1)=Z$  (resp. $\on(\hat Z_2)=Z$).  Alors $f_*(\hat Z_1)=\hat Z_2$, et si $p \in \mbox{P}$, $(R_p)_*(\hat Z)=\hat {Z_p}$ , o\`u $Z_p:={ (\Ad p^{-1}).Z}$.  On en d\'eduit la propri\'et\'e d'\'equivariance tr\`es importante suivante.  Pour tout $\xi \in {\mathcal W}_{\hx}$, et $p \in \mbox{P}$, on a $(\Ad p).\xi \in {\mathcal W}_{f(\hx).p^{-1}}$ et~:
 \begin{equation}
 \label{equ.equivariance-exponentielle}
  f(\exp(\hx,\xi)).p^{-1}=\exp(f(\hx).p^{-1}, (\Ad p).\xi).
  \end{equation}   

%\subsection{Les structures conformes de dimension $\geq 3$ vues comme structures de Cartan}

\subsection{Sous-vari\'et\'es totalement g\'eod\'esiques conformes}
\label{sec.geodesiques}
Dans toute cet\-te section, $(M,[g])$ d\'esigne une structure conforme pseudo-riemannienne de dimension $n \geq 3$, et de signature $(p,q)$.  On interpr\`ete cette structure conforme comme une g\'eom\'etrie de Cartan model\'ee sur $\Ein^{p,q}$, et l'on appelle $(M,\hm,\omega)$ le fibr\'e normal de  Cartan.  % canonique d\'efini par $(M,g)$.   L'int\'er\^et de ce point de vue est qu'il donne  une notion de {\it d\'eveloppement} de toute courbe 

La connexion de Cartan $\omega$  d\'efinit une notion de transport parall\`ele sur le fibr\'e $\hm$~: si $\alpha:  [0,1] \to \hm$ est une courbe lisse et si $\xi \in \lieg$ est un vecteur, alors $\omega_{\alpha(t)}^{-1}(\xi)$ d\'efinit un champ de vecteur le long de la courbe $\alpha$, que l'on qualifiera de ``parall\`ele".  Les courbes de $\hm$ dont le vecteur tangent est parall\`ele sont les $t \mapsto \exp(\hx,\xi)$, $\hx \in \hm$, $\xi \in \lieg$.  On d\'efinit  les {\it segments g\'eod\'esiques conformes}  (param\'etr\'es) de $M$ comme les projections sur $M$ de courbes de la forme $t \mapsto \exp(\hx, \xi)$, $\hx \in \hm$, $\xi \in \lien^-$ (noter que l'on se limite aux vecteurs $\xi$ horizontaux).  Par exemple, sur la sph\`ere standard riemannienne ${\bf S}^n$, les segments g\'eod\'esiques conformes sont les arcs de grands cercles, ainsi que leurs images par les transformations de M\"obius.
 Nous ne d\'etaillerons pas les propri\'et\'es de ces courbes dans cet article (voir \cite{cap}, \cite{ferrand3}).

De mani\`ere g\'en\'erale, \'etant donn\'e  ${\mathcal S} \subset \lieg$ un sous-espace vectoriel, on peut se demander si la distribution de $T\hm$ d\'efinie par $\omega^{-1}({\mathcal S})$ admet une ou plusieurs feuilles int\'egrales. Autrement dit, existe-t-il une sous-vari\'et\'e $\hat S \subset \hm$ telle que $\omega(T\hat S)=\{ {\mathcal S} \}$; on dit alors que $\hat S$ est une sous-vari\'et\'e {\it parall\`ele} de $\hm$.  Notons que n\'ecessairement, si $\hat S$ est une sous-vari\'et\'e int\'egrale de $\omega^{-1}(\mathcal S)$ passant par $\hx$, alors pour ${\mathcal U}$ un voisinage de $0$ suffisamment petit, $\exp(\hx,{\mathcal S} \cap {\mathcal U}) \subset \hat S$.

L'existence de  sous-vari\'et\'es parall\`eles de dimension $>1$ traduit  souvent des propri\'et\'es g\'eom\'etriques particuli\`eres de la structure conforme $(M,[g])$.  \`A titre d'exemple, citons le~:

\begin{lemme}
\label{lem.horizontal.integrable.plat}
Si la distribution de $T \hm$ d\'efinie par $\omega^{-1}(\lien^-)$ admet une feuille int\'egrale $\hat S$ passant par $\hx \in \hm$, alors il existe un ouvert de $M$ contenant $x:=\pi_M(\hx)$ qui est localement conform\'ement plat.
\end{lemme}

\begin{preuve} soit $Z_1, \ldots, Z_n$ une base de $\lien^-$, et $\hat Z_1, \ldots, \hat Z_n$ les champs $\omega$-constants de $\hm$ associ\'es, alors on a en chaque point $\hy$ de $\hat S$ la relation~:
$$ K(\hat Z_i, \hat Z_j)=[Z_i,Z_j] - \omega([\hat Z_i, \hat Z_j])=- \omega([\hat Z_i, \hat Z_j]).$$
Comme $\omega(T\hat S) \subset \lien^-$ par d\'efinition, on a qu'en chaque $\hy \in \hat S$,  $K(\hat Z_i, \hat Z_j) \in \lien^-$.  La connexion de Cartan \'etant normale, sa courbure est \`a valeurs dans $\R \oplus \oo(p,q) \oplus \lien^+$.  On conclut que $K=0$ sur $\hat S$, et comme $\hat S$ se projette sur un ouvert contenant $\pi_M(\hx)$, le lemme s'ensuit. \end{preuve}

Par analogie avec les g\'eod\'esiques conformes, on va maintenant d\'efinir les sous-vari\'et\'es totalement g\'eod\'esiques conformes de $M$ comme des projections sur $M$ de certaines sous-vari\'et\'es parall\`eles de $\hm$.

\begin{definition}[Sous-vari\'et\'es totalement g\'eod\'esiques conformes]
Soit $\Sigma \subset M$ une sous-vari\'et\'e.  On dit que $\Sigma$ est totalement g\'eod\'esique conforme lorsque~:
\begin{itemize}
\item{Il existe une sous-alg\`ebre de Lie $\lieh:=\lien_1 \oplus \liep_1$ dans $ \lieg$, avec $\lien_1 \subset \lien^-$,  $\liep_1 \subset \liep$, telle que $\omega^{-1}(\lieh)$ admette une feuille int\'egrale $\hat S$ dans $\hm$.}
\item{La sous-vari\'et\'e $\Sigma$ s'\'ecrit $\Sigma=\pi_M(\hat S)$.}
\end{itemize}
\end{definition}

Comme nous l'avons dit, il n'y a  g\'en\'eriquement pas de sous-vrai\'et\'e totalement g\'eod\'esique de dimension $\geq 2$.

 \subsection{Existence locale de m\'etriques Ricci-plates dans la classe conforme}

Soit $(M,\hm,\omega)$ le fibr\'e normal de Cartan d'une vari\'et\'e pseudo-riemannienne conforme $(M,[g])$.  Nous avons vu au Lemme \ref{lem.horizontal.integrable.plat} que si $\lieh$ est une sous-alg\`ebre de Lie de $\lieg$, l'existence d'une feuille int\'egrale \`a la distribution $\omega^{-1}(\lieh)$ sur $T\hm$ pouvait avoir des cons\'equences g\'eom\'etriques int\'eressantes.   Un  autre cas remarquable est donn\'e par la proposition suivante~:

\begin{proposition}
\label{prop.ricci.plate}
Si la distribution $\omega^{-1}(\lien^- \oplus \lies)$ admet une feuille int\'egrale $\hat M_0 \subset \hm$ passant par $\hx \in \hm$, alors il existe un voisinage $U$ de $x=\pi_M(\hx)$, et une m\'etrique Ricci-plate dans la classe conforme $[g]_{|U}$.
\end{proposition}  

\begin{preuve} dans la preuve, on va identifier $\lieg/\liep$ \`a $\lien^-$, via un isomorphisme qui commute \`a l'action adjointe de $\mbox{Z} \times \mbox{S}$.  On  consid\'erera alors $\iota_{\hx}$ comme un isomorphisme entre $\lien^-$ et $T_xM$.  Apr\`es cette identification, la relation~:
\begin{equation}
\label{equn.relation}
 \iota_{\hx.p^{-1}}((\Ad p).{\xi})=\iota_{\hx}(\xi)
\end{equation}  
est encore valable  lorsque $\xi \in \lien^-$, et $p \in \mbox{Z} \times \mbox{S} \subset \mbox{P}$.  

La feuille int\'egrale $\hat M_0$, si elle est choisie maximale, est stable par l'action \`a droite de $\mbox{S}^o$ sur $\hm$.  Quitte \`a prendre le satur\'e par l'action de $ \mbox{S}$, nous supposerons dans ce qui suit que $\hm_0$ est stable par l'action de $ \mbox{S}$. Soit ${\mathcal U}$  un voisinage de $0$ dans $\lieg$, assez petit pour que $\xi \mapsto \pi_M(\exp(\hx,\xi))$ soit un diff\'eomorphisme de ${\mathcal U} \cap \lien^-$ sur son image.  On note $\hat \Sigma:=\exp(\hx,\lien^- \cap {\mathcal U})$, et $\hat U$ le satur\'e de $\hat \Sigma$ par l'action de $ \mbox{S}$; c'est un ouvert de $\hat M_0$ diff\'eomorphe au produit $\hat \Sigma \times  \mbox{S}$.  On pose \'egalement $U:=\pi_M(\hat \Sigma)$, qui est un voisinage de $x:=\pi_M(\hx)$.   Les fibres de la restriction $\overline \pi_M$ de $\pi_M$ \`a $\hat U$ sont exactement les orbites de $ \mbox{S}$, et $\overline \pi_M : \hat U \to U$ est un $ \mbox{S}$-fibr\'e principal (ou de mani\`ere \'equivalente, un $\mbox{O}(p,q)$-fibr\'e principal).

La donn\'ee de la sous-vari\'et\'e $\hat U$ permet de d\'efinir une m\'etrique $h$ sur $U$ comme suit~: $\text{pour tout } y \in U \text{ et } u,v \in T_yM, \ h_y(u,v):=\lambda^{p,q}(\iota_{\hy}^{-1}(u),\iota_{\hy}^{-1}(v))$.  V\'erifions que $h$ est  d\'efinie sans ambig\"uit\'e; si $\hy^{\prime} \in \hat U$ est dans la m\^eme fibre que $\hy$, alors $\hy^{\prime}=\hy.p$ pour un \'el\'ement $p \in  \mbox{S}$.   Alors par la relation (\ref{equn.relation}), pour tout $u \in T_yM$, $\iota_{\hy}^{-1}(u)=(\Ad p).\iota_{\hx.p}^{-1}(u)$.  Comme le produit $\lambda^{p,q}$ est invariant par l'action adjointe de $ \mbox{S}$, on a bien $\lambda^{p,q}(\iota_{\hy}^{-1}(u),\iota_{\hy}^{-1}(v))=\lambda^{p,q}(\iota_{\hy^{\prime}}^{-1}(u),\iota_{\hy^{\prime}}^{-1}(v))$.

Soit $\hat R$ le fibr\'e des rep\`eres orthonorm\'es de $h$.  Si $(\xi_1, \ldots , \xi_n)$ est une base orthonorm\'ee de $\lien^-$, on d\'efinit un isomorphisme de fibr\'es $\varphi$ entre $\hat U$ et $\hat R$ de la fa\c{c}on suivante~:  pour $\hy \in \hat U$, on pose $\varphi(\hy):=(\iota_{\hy}(\xi_1), \ldots , \iota_{\hy}(\xi_n))$. 

Appelons maintenant $\tilde \omega$ la restriction de la connexion de Cartan $\omega$ \`a $T \hat U$.  Il s'agit d'une connexion de Cartan sur le $\mbox{S}$-fibr\'e $\hat U$, \`a valeurs dans $\lien^- \oplus \lies$, qui d\'efinit la m\'etrique $h$.  Le point essentiel est que si $\hat Z_i$ et $\hat Z_j$ sont deux champs de vecteurs $\omega$-constants, avec $\omega(\hat Z_i) \in \lien^-$ et $\omega(\hat Z_j) \in \lien^-$, alors d\`es que $\hy \in \hat U$, les vecteurs $\hat Z_i(\hy), \hat Z_j(\hy)$ et $[\hat Z_i,\hat Z_j](\hy)$ appartiennent tous trois \`a $ T_{\hy}\hat U$.  On en d\'eduit que les fonctions courbures $\kappa : \hm \to \Hom (\Lambda^2(\lien^-),\lies \oplus \lien^+)$ et $\tilde \kappa: \hm \to \Hom (\Lambda^2(\lien^-),\lien^- \oplus \lies)$ de $\tilde \omega$, {\it  co\"{\i}ncident sur $\hat U$}.  Ainsi, si $\kappa_{\lies}$ est la composante selon $\lies$ de $\kappa$, on a $\kappa=\kappa_{\lies}=\tilde \kappa$ en restriction \`a $\hat S$.  La premi\`ere cons\'equence est que $\tilde \kappa$ est en fait \`a valeurs dans $\Hom (\Lambda^2(\lien^-),\lies)$, ce qui signifie qu'elle est sans torsion.  C'est donc que $\tilde \omega$ est la connexion de Levi-Civita de $h$ (voir \cite{sharpe} chap 6. pour une pr\'esentation du ``probl\`eme" d'\'equivalence dans le cadre m\'etrique).  Maintenant, le fait que $\kappa_{\lies}=\tilde \kappa$ sur $\hat U$ signifie que sur $U$, le tenseur de courbure de $h$ co\"{\i}ncide avec  le tenseur de Weyl.  En effet, le premier est donn\'e par $R_y(u,v)w=[\tilde \kappa_{\hy}(\iota_{\hy}^{-1}(u),\iota_{\hy}^{-1}(v)),\iota_{\hy}^{-1}(w)]$, tandis que le second s'exprime comme $W_y(u,v)w=[\kappa_{\lier}(\iota_{\hy}^{-1}(u),\iota_{\hy}^{-1}(v)),\iota_{\hy}^{-1}(w)]$.  La courbure de Ricci de $h$ doit alors s'annuler (voir \cite{besse} Th\'eor\`eme 1.114, p. 47, et \cite{sharpe} Proposition 1.4, p. 229 pour la d\'ecomposition de l'espace des tenseurs de courbure en composantes irr\'eductibles sous l'action de $\mbox{O}(p,q)$). \end{preuve}

\section{Familles normales d'immersions conformes}
\label{sec.dynamique-holonomie}
Dans toute cette section, $(M,[g])$ et $(N,[h])$ d\'esignent deux  structures conformes pseudo-riemanniennes de m\^eme signature $(p,q)$, avec $p+q \geq 3$. Jusqu'\`a la fin de l'article, on d\'esignera par $(M,\hm,\om)$ et $(N,\hn,\on)$ les fibr\'es normaux de Cartan associ\'es aux structures conformes sur $M$ et $N$ respectivement. 

\subsection{Holonomie d'une suite d'immersions conformes}
\label{sec.holonomie-suite}
\ \\
Soit  $f_k : (M,g) \to (N,h)$, $k \in \N$, une famille d'immersions conformes.  On suppose que  $x \in M$  est tel que $f_k(x)$ soit relativement compacte  dans $N$.  On dit alors qu'une suite  $(h_k)$ de $\mbox{P}$ est {\it une suite d'holonomie de $(f_k)$ en $x$} s'il existe une suite $(\hx_k)$ dans la fibre de $x$, contenue dans un compact de $\hm$, et  telle que $f_k(\hx_k).h_k^{-1}$ soit \'egalement contenue dans un compact de $\hn$. Remarquons que, sous l'hypoth\`ese o\`u $f_k(x)$ est relativement compacte dans $M$, il existe toujours au moins une suite d'holonomie associ\'ee \`a $(f_k)$.  

\subsubsection{Suites \'equivalentes}
\label{sec.suites.equivalentes}
La notion de suite d'holonomie est stable ``par perturbation compacte" : si $(h_k)$ est une suite d'holonomie de $(f_k)$ en $x$, alors il en va de m\^eme pour toute suite $h_k^{\prime}=l_1(k)h_kl_2(k)$, o\`u $l_1(k)$ et $l_2(k)$ sont des suites  relativement compactes de $\mbox{P}$.  On dit alors que $(h_k)$ et $(h_k^{\prime})$ {\it sont \'equivalentes}.  
Du fait que l'action de $\mbox{P}$ sur $\hm$ et $\hn$ est propre, il est facile de v\'erifier que r\'eciproquement, deux suites d'holonomies de $(f_k)$ en $x$ sont toujours \'equivalentes.
%Ce sont les propri\'et\'es de la classe d'\'equivalence d'une suite d'holonomie qui vont \^etre utiles
Donc, ce qui a vraiment un sens, c'est la classe d'\'equivalence des suites d'holonomies de $(f_k)$ en un point $x$. Mais dans tout l'article, et par abus de langage, on dira fr\'equemment~: {\it soit $(h_k)$ l'holonomie de $(f_k)$ en $x$}.  On entendra par l\`a que $(h_k)$ est un repr\'esentant de la classe d'\'equivalence de suites d'holonomies en $x$. On sera souvent amen\'e \`a  changer une suite d'holonomie en une autre  qui lui est \'equivalente.

\subsection{Notion de stabilit\'e.}
Soit  $f_k : (M,g) \to (N,h)$ 
une suite d'immersions  conformes. %Soit $x$ un point de $M$ tel que la suite $f_k(x)$ admette une limite $z_{\infty}$ dans $N$.  
\begin{definition}[Stabilit\'e]
\label{defi.stabilite}
On dit que la  suite $(f_k)$ est stable en $x \in M$ si $f_k(x)$ converge vers une limite $z_{\infty} \in N$, et si pour toute suite $(x_k)$ de $M$ convergeant vers $x$, $f_k(x_k)$ tend \'egalement vers $z_{\infty} $.  La suite $(f_k)$ est dite fortement stable en $x$ s'il existe un voisinage $U$ de $x$ dans $M$ tel que $f_k(\overline{U})$ converge vers  $z_{\infty} \in N$ pour la topologie de Hausdorff.
\end{definition}

On va \'egalement d\'egager une notion de stabilit\'e pour les suites de $\mbox{P}$ :
\begin{definition}
\label{defi.opq-stabilite} 
Une suite $(h_k)$ de $\mbox{P}$ est dite stable si c'est une suite de $\mbox{A}^+$. De mani\`ere \'equivalente, $(h_k)$ est stable si elle s'\'ecrit~:
$$h_k=\text{\diag}(\lambda_1(k), \ldots , \lambda_n(k)) \in \R_+^* \times \mbox{O}(p,q),$$  avec $\lambda_1(k) \geq \ldots \geq  \lambda_n(k) \geq 1$; en particulier  les suites $(\frac{1}{\lambda_i(k)})$, $i=1, \ldots, n$, sont born\'ees. La suite $(h_k)$ est dite fortement stable lorsqu'elle est stable et que de plus les suites  $(\frac{1}{\lambda_i(k)})_{k \in \N}$, $i=1, \ldots, n$ tendent vers $0$.  
%\end{enumerate}
\end{definition}

%\'Etant donn\'ee une suite $(h_k)$ de $\mbox{P}$, on peut consid\'erer $(h_k)$ comme une suite d'immersions conformes de $\Ein^{p,q}$, et se demander si elle est stable en $o$ au sens de la d\'efinition \ref{defi.stabilite}.  On peut \'egalement s'interroger sur sa stabilit\'e au sens de la d\'efinition \ref{defi.opq-stabilite}.  En fait, ces deux notions co\"{\i}ncident, comme le montre le lemme qui suit.  

Le lemme suivant va montrer que la propri\'et\'e, pour une suite $(f_k)$ comme ci-dessus, d'\^etre stable en $x$, se lit sur l'holonomie $(h_k)$ de $(f_k)$ en $x$.  En particulier, une suite $(h_k)$ de $\mbox{P}$, vue comme suite d'immersions conformes de $\Ein^{p,q}$, est stable en $o$ (au sens de la d\'efinition \ref{defi.stabilite}) si et seulement si elle est stable au sens de \ref{defi.opq-stabilite}.  Il n'y a donc pas d'ambig\"uit\'e dans la terminologie. 

\begin{lemme}
\label{lem.holonomie-stable}
La suite $(f_k)$ est stable en $x \in M$ (resp. fortement stable) si et seulement si $f_k(x)$ converge vers $z_{\infty} \in N$ et s'il existe en $x$ une suite d'holonomie $(h_k)$ qui soit  stable (resp.  fortement stable).  
\end{lemme}

\begin{preuve} on commence par supposer que $f_k(x)$ converge vers $z_{\infty}$ et que $(f_k)$ est stable en $x$.  On se donne $(h_k)$ une suite d'holonomie de $(f_k)$ en $x$, et l'on exprime $h_k$ sous forme affine :
$$ h_k=\sigma_kL_k+T_k$$
o\`u $(\sigma_k), (L_k)$ et $(T_k)$ sont  des suites de $\R_+^*, \mbox{O}(p,q)$ et $\R^{p,q}$ respectivement.  On peut \'ecrire $L_k=L_{1k}D_kL_{2k}$ o\`u $D_k$ est un \'el\'ement de  $\mbox{O}(p,q)$, de la forme $D_k=\diag (\lambda_1(k), \ldots, \lambda_n(k))$, avec $\lambda_1(k) \geq \ldots \geq \lambda_n(k) >0$.  Les suites   $(L_{1k})$, $(L_{2k})$ appartiennent \`a  $\mbox{O}(p,q)$ et sont relativement compactes.  Ainsi, quitte \`a remplacer $(h_k)$ par une suite \'equivalente, on peut supposer que $(h_k)$ est \'egale \`a $\sigma_k D_k(Id+\tau_k)$.    Notre but dans un premier temps, est de montrer~:
\begin{fait}
\label{fait.trans.borne}
 Si $(f_k)$ est stable en $x$, alors la suite $(\tau_k)$ doit \^etre born\'ee.
\end{fait}

\begin{preuve} d'apr\`es les expressions matricielles (\ref{eq.matricielle}) et (\ref{eq.matricielle2}) donn\'ees en section \ref{sec.matriciel}, la suite $(h_k)$ s'exprime dans $\mbox{O}(p+1,q+1)$ sous la forme :

$$ \left( \begin{array}{ccc}
\sigma_k& -\sigma_k\tau_k^t.J_{p,q}   &  -\frac{\sigma_k}{2} Q^{p,q}(\tau_k)\\
  & D_k  & D_k.\tau_k \\
  &   & \sigma_k^{-1}
\end{array} \right).$$

Supposons par l'absurde que  $(\tau_k)$ ne soit pas born\'ee. On \'ecrit :
$$\tau_k^t:=(\tau_1(k), \ldots , \tau_n(k)).$$
 Quitte \`a consid\'erer une suite extraite de $(f_k)$, il existe $i \in \{ 1, \ldots,  n \}$ tel que $|\tau_i(k)| \to \infty$.  Soit $j=n+1-i$ si $i \in \{ 1, \ldots , p \} \cup \{ q+1, \ldots ,n \}$, et $j=i$ sinon.  Notons $\Delta$ le projectivis\'e de $\text{Vect}(e_0,e_j)$ sur $\Ein^{p,q}$.   Alors $h_k$ pr\'eserve $\Delta$ et son action se fait via la transformation de $PSL(2,\R)$  suivante :
$$ \left( \begin{array}{cc} \sigma_k & -\sigma_k \tau_i(k)\\
0 & \lambda_j(k)\end{array} \right).$$ 
Ainsi, il existe un param\'etrage projectif de $\Delta$ au voisinage de $o$, donn\'e par $s \mapsto \pi_G(e^{s \zeta})$, tel que $e^{(\Ad h_k).s\zeta}=e^{\alpha_k(s)\zeta}.p_k(s)$, o\`u $\alpha_k(s)=\frac{\lambda_j(k)s}{\sigma_k-\sigma_k\tau_i(k)s}$ est d\'efinie sur $I_k:= \R \setminus \{  \frac{1}{\tau_i(k)}\}$, et $p_k : I_k \to P$ est une courbe valant $1_G$ en $0$.  
Comme $(h_k)$ est une suite d'holonomie de $(f_k)$, quitte \`a consid\'erer une suite extraite de $(f_k)$, on peut supposer qu'il existe une suite $(\hx_k)$ qui converge vers $\hx$, telle que $f_k(\hx_k).h_k^{-1}$ converge vers $\hz \in \hn$. On choisit $\lambda >0$ assez petit pour  que $\hy \mapsto \exp(\hy, \lambda \zeta)$ soit bien d\'efini sur un voisinage de $\hz$, et que de plus $\pi_N(\exp(\hz,\lambda \zeta)) \not = z_{\infty}$  (c'est possible puisque $\zeta$ est transverse \`a $\liep$).  

Quitte \`a consid\'erer \`a nouveau une suite extraite de $(f_k)$, on peut supposer que $\frac{1}{c\tau_i(k)}$ est de signe constant, par exemple positif.  L'image $\alpha_k(]0,\frac{1}{\tau_i(k)}[)$ est $]- \infty, 0[$ ou $]0, \infty[$.  L\`a encore, quitte \`a extraire une sous-suite, on peut supposer que c'est toujours le m\^eme intervalle, par exemple $]0, \infty[$  (les autres cas se tra\^{\i}tent  de fa\c{c}on similaire).  Ainsi, pour tout $k$, il existe $s_k \in \ ]0,\frac{1}{\tau_i(k)}[$ tel que $\alpha_k(s_k)=\lambda$.  Par ailleurs, $s_k \to 0$ puisque $|\tau_i(k)| \to \infty$.

 Pour $k$ suffisamment grand, on peut \'ecrire :

\begin{equation}
\label{equn2}
f_k(\exp(\hx_k,s_k \zeta)).h_k^{-1}=\exp(f_k(\hx_k).h_k^{-1},\lambda \zeta).p_k(s_k).
 \end{equation}

Expliquons pourquoi cette relation est vraie. Soit $(S,\hat S,\omega^S)$  le fibr\'e normal de Cartan d'une vari\'et\'e pseudo-riemannienne conforme $S$ de signature $(p,q)$. Si $\lambda : I \to \hat S$ et $p : I \to P$  sont deux  courbes de classe $C^1$, et si l'on pose $\gamma(t)=\lambda(t).p(t)$, alors (voir  {\cite[p. 208]{sharpe}})~:
\begin{equation}
\label{equ.formule-Cartan}
\omega^S(\gamma^{\prime}(t))=(\Ad p(t))^{-1}.\omega^S(\lambda^{\prime}(t))+\omega_G(p^{\prime}(t)),
\end{equation}
o\`u $\omega_G$ d\'esigne la forme de Maurer-Cartan sur le groupe $G=\mbox{O}(p+1,q+1)$.
On peut appliquer  l'\'equation (\ref{equ.formule-Cartan}), dans le mod\`ele $(G/P,G,\omega_G)$, \`a l'\'egalit\'e $e^{(\Ad h_k).s\zeta}=e^{\alpha_k(s)\zeta}.p_k(s)$, et l'on obtient~:
$$(\Ad h_k).\zeta=\alpha_k^{\prime}(s)(\Ad p_k(s))^{-1}.\zeta + \omega_G(p_k^{\prime}(s)).$$
Posons :
$$\gamma_1(s):=f_k(\exp(\hx_k,s \zeta)).h_k^{-1}=\exp(f_k(\hx_k).h_k^{-1},s(\Ad h_k).\zeta)$$
et
$$\gamma_2(s):=\exp(f_k(\hx_k).h_k^{-1},\alpha_k(s) \zeta).p_k(s).$$
On a
$$ \on(\gamma_1^{\prime}(s))=(\Ad h_k).\zeta$$
et l'\'equation (\ref{equ.formule-Cartan}) 
fournit~:
$$ \on(\gamma_2^{\prime}(s))= \alpha_k^{\prime}(s)(\Ad p_k(s))^{-1}.\zeta + \omega_G(p_k^{\prime}(s)).$$

Les courbes $s \mapsto \gamma_1(s)$ et $s \mapsto \gamma_2(s)$ satisfont donc la m\^eme \'equation diff\'erentielle du premier ordre, et valent toutes deux $f_k(\hx_k).h_k^{-1}$ en $0$, donc elles sont \'egales, ce qui justifie la relation (\ref{equn2}).
%Ceci est une cons\'equence du lemme suivant (voir \cite{frances-melnick}, lemme .....).

En projetant l'\'egalit\'e (\ref{equn2}) sur $M$ et $N$, on obtient~:
$$ f_k(\pi_M(\exp(\hx_k,s_k \zeta)))=\pi_N(\exp(f_k(\hx_k).h_k^{-1},\lambda \zeta)).$$

Or   $\pi_M(\exp(\hx_k,s_k \zeta))$ tend vers $x$ puisque $s_k \to 0$.\\
 En revanche,  $\pi_N(\exp(f_k(\hx_k).h_k^{-1},\lambda \zeta))$  tend vers $\pi_N(\exp(\hz,\lambda \zeta)) \not = z_{\infty}$.  Cela contredit le fait que $(f_k)$ est stable en $x$. \end{preuve}

On d\'eduit de ce qui pr\'ec\`ede que $(\tau_k)$ est une suite born\'ee, et que par cons\'equent, $(h_k)$ est \'equivalente dans $\mbox{P}$ \`a une suite de la forme $h_k=\text{diag}(\lambda_1(k), \ldots , \lambda_n(k)) \in \R_+^* \times \mbox{O}(p,q)$, o\`u $\lambda_1(k) \geq \ldots  \geq  \lambda_n(k) >0$. On veut \`a pr\'esent montrer~:
\begin{fait}
\label{fait.dist.borne}
Si $(f_k)$ est stable (resp. fortement stable) en $x$, alors les suites $\frac{1}{\lambda_i(k)}$ sont born\'ees (resp. tendent vers $0$).  
\end{fait}

\begin{preuve} supposons par l'absurde que le Fait \ref{fait.dist.borne} n'a pas lieu. Il existe  alors une suite $(\xi_m)$ de vecteurs de $\lien^-$ qui tend vers $0$, et une suite extraite $(h_{k_m})$ telle que $(\Ad h_{k_m}).\xi_m$ tende vers $\xi \not = 0$.  Quitte \`a multiplier $\xi_m$ par un $\epsilon>0$ assez petit, on peut supposer que l'application exponentielle est bien d\'efinie sur un voisinage de $(\hz,\xi)$ dans $\hn \times \lieg$, et que $\exp(\hz,\xi) \not = z_{\infty}$.  On \'ecrit alors~:
$$ f_{k_m}(\pi_M(\exp(\hx_{k_m},\xi_m)))=\pi_N(\exp(f_{k_m}(\hx_{k_m}).h_{k_m}^{-1}, (\Ad h_{k_m}). \xi_m)).$$
Or $\pi_M(\exp(\hx_{k_m},\xi_m)$ tend vers $x$.\\
En revanche,  $\pi_N(\exp(f_{k_m}(\hx_{k_m}).h_{k_m}^{-1}, (\Ad h_{k_m}).\xi_m))$ tend vers $\pi_N(\exp(\hz,\xi) \not = z_{\infty})$~: contradiction avec la stabilit\'e de $(f_k)$ en $x$.

Maintenant, si l'une des suites $\frac{1}{\lambda_i(k)}$ ne tend pas vers $0$, alors il existe une suite extraite $\frac{1}{\lambda_i(k_m)}$ ayant une limite non nulle $\lambda$.  On choisit $\xi_i \in \lien^-$ un vecteur propre de $\Ad h_k$ pour la valeur propre $\frac{1}{\lambda_i(k)}$.  Notons que $\Ad h_k$ est diagonale dans une base ind\'ependante de $k$, donc $\xi_i$ peut \^ etre choisi ind\'ependant de $k$.  Ainsi :
$$ f_{k_m}(\pi_M(\exp(\hx_{k_m},s \xi_i)))=\pi_N(\exp(f_{k_m}(\hx_{k_m}).h_{k_m}^{-1},\frac{1}{\lambda_i(k_m)} s \xi_i)).$$
Pour $s$ petit, $\pi_M(\exp(\hx_{k_m},s \xi_i))$ converge vers un point proche de $x$, tandis que $\pi_N(\exp(f_{k_m}(\hx_{k_m}).h_{k_m}^{-1},\frac{1}{\lambda_i(k_m)} s \xi_i))$ converge vers $\exp_N(\hz,\lambda s \xi_i)$, qui est diff\'erent de $z_{\infty}$.  

La suite $(f_k)$ n'est pas fortement stable en $x$ dans ce cas. \end{preuve}

Le Fait \ref{fait.dist.borne} affirme que  si $(f_k)$ est stable en $x$ (resp, fortement stable), alors $(h_k)$ est une suite stable (resp. fortement stable) de $\mbox{P}$.
Il nous reste \`a  montrer que r\'eciproquement, si $f_k(x)$ tend vers $z_{\infty}$ et si $(f_k)$ admet une suite d'holonomie $(h_k)$ en $x$ qui est une suite stable (resp.  fortement stable) de $\mbox{P}$, alors $(f_k)$ est stable (resp.  fortement stable) en $x$.  Supposons pour commencer que  la suite d'holonomie $(h_k)$ est stable, et s'\'ecrit $\text{diag} (\lambda_1(k), \ldots , \lambda_n(k))$, o\`u les suites $\frac{1}{\lambda_i(k)}$ sont born\'ees. On consid\`ere une suite $(\hx_k)$ de la fibre de $x$ contenue dans un compact de $\hm$, de sorte que $f_k(\hx_k).h_k^{-1}$ soit  contenue dans un compact de  $\hn$.  Soit $(y_k)$ une suite qui tend vers $x$.  Alors pour $k$ assez grand, on peut \'ecrire $y_k=\pi_M(\exp(\hx_k,\xi_k))$, pour une certaine suite $(\xi_k)$ de $\lien^-$ qui tend vers $0$.  On obtient que :
$$ f_k(\exp(\hx_k,\xi_k)).h_k^{-1}=\exp(f_k(\hx_k).h_k^{-1},(\Ad h_k).\xi_k).$$
Or les $\frac{1}{\lambda_i(k)}$ \'etant born\'ees, on a $(\Ad h_k).\xi_k \to 0$. En projetant sur la vari\'et\'e $N$, on obtient bien $f_k(y_k) \to z_{\infty}$ : la suite $(f_k)$ est stable en $x$.

Si de plus, on sait que $\lim_{k \to \infty}\frac{1}{\lambda_i(k)} = 0$ pour  $i=1, \ldots ,n$, alors pour tout ensemble relativement compact $K$ de $\lien^-$, on a $(\Ad h_k).K \to   0 $  (la limite \'etant prise pour la topologie de Hausdorff).  On choisit \`a pr\'esent $U$ un voisinage suffisamment petit de $x$, de sorte que pour tout $k$ assez grand, il existe ${\mathcal U}_k \subset \lien^-$ un voisinage de $0$ tel que  $\xi \mapsto \pi_M( \exp(\hx_k,\xi))$ soit un diff\'eomorphisme de ${\mathcal U}_k$ sur $U$.  Par relative compacit\'e de la suite $(\hx_k)$, si $U$ est pris assez petit, les ${\mathcal U}_k$ sont tous inclus dans un compact de $\lien^-$.  On obtient alors :
$$ f_k(\exp(\hx_k,{\mathcal U}_k)).h_k^{-1}=\exp(f_k(\hx_k).h_k^{-1},(\Ad h_k).{\mathcal U}_k).$$
En projetant cette relation sur $M$ et $N $, on aboutit \`a $\lim_{k \to \infty}f_k(U)=z_{\infty}$, la limite \'etant prise pour la topologie de Hausdorff.  La suite $(f_k)$ est bien fortement stable en $x$. \end{preuve}

\subsection{Caract\'erisation des familles normales~: preuve du th\'eor\`eme \ref{thm.normal}}
\label{sec.caracterisation}
Ici, $(M,[g])$ et $(N,[h])$ sont deux structures pseudo-riemanniennes de signature $(p,q)$, $p+q=3$, et ${\mathcal F}$ est une famille d'immersions conformes de $(M,[g])$ dans $(N,[h])$. Il est clair que s'il  existe un voisinage $U$ de $x$ tel que la famille ${\mathcal F}_{|U}$ est d'adh\'erence compacte dans $\Cont (U,N)$, alors d'une part, $E=\{f(x) \ | \ f \in {\mathcal F}  \}$ est relativement compact dans $N$, et de plus ${\mathcal F}$ est \'equicontinue en $x$.

Nous allons   supposer r\'eciproquement  qu'il existe $x \in M$ tel que $E=\{f(x) \ | \ f \in {\mathcal F}  \}$ soit relativement compact dans $N$, et que de plus ${\mathcal F}$ soit \'equicontinue en $x$.  On veut montrer l'existence de $U \subset M$ un ouvert contenant $x$, tel que  toute suite de ${\mathcal F}_{|U}$ admet une sous-suite qui converge vers une application lisse, la convergence \'etant  $C^{\infty}$ sur les compacts de $U$.
On commence par choisir un compact ${\mathcal K} \subset \hn$ qui se projette surjectivement sur l'adh\'erence 
$\overline{E}$ de $E$ dans $N$.  On choisit \'egalement $\hx \in \hm$ dans la fibre de $x$.  Pour tout $f \in {\mathcal F}$, on choisit  $p(f) \in \mbox{P}$ tel que $f(\hx).p(f)^{-1} \in {\mathcal K}$.  On h\'erite ainsi d'une application $p:{\mathcal F} \to P$.  Nous allons voir que l'hypoth\`ese d'\'equicontinuit\'e de ${\mathcal F}$ en $x$ va avoir des cons\'equences sur l'image de $\mbox{P}$.

Tout d'abord, on \'ecrit $p(f)$ sous forme affine $p(f)=\sigma(f)L(f)+T(f)$, o\`u $\sigma(f) \in \R_+^*$, $L(f) \in \mbox{O}(p,q)$, et $T(f) \in \R^{p,q}$.  On note $K_0$ le compact maximal de $\mbox{O}(p,q)$.  En effectuant une d\'ecomposition de Cartan de $L(f)$, on \'ecrit  $p(f)=L_1(f)A(f)L_2(f)+T(f)$, o\`u $L_1(f)$ et $L_2(f)$ sont dans $K_0$ et $A(f)= \diag (\lambda_1(f), \ldots ,\lambda_n(f)) \in  \R_+^* \times \mbox{O}(p,q)$ v\'erifie $\lambda_1(f) \geq  \ldots \lambda_n(f)>0$.  On \'ecrit finalement $p(f)=L_1(f)A(f)(id + \tau(f))L_2(f).$  On commence alors par remarquer qu'il existe un compact $K_1 \subset \mbox{P}$ tel que l'application $f \mapsto id + \tau(f)$ soit \`a valeurs dans $K_1$.  Si tel n'\'etait pas le cas, il existerait une suite $(f_k)$ de ${\mathcal F}$ telle que $\tau_k:=\tau(f_k)$ ne soit pas born\'ee dans $\R^{p,q}$.  Quitte \`a consid\'erer une suite extraite de $(f_k)$, on peut supposer que $f_k(x)$ tend vers $z_{\infty} \in N$.  Mais nous avons vu au Fait \ref{fait.trans.borne} que dans ce cas, $(f_k)$ n'\'etait pas stable en $x$.  Ceci contredirait l'\'equicontinuit\'e de la famille ${\mathcal F}$ en $x$.  De la m\^eme mani\`ere, il existe un compact $K_2 \subset \mbox{M}(n,\R)$ de sorte que l'application $f \mapsto A(f)$ soit \`a valeur dans $K_2$. Si tel n'\'etait pas le cas, on pourrait trouver une suite $(f_k)$ de ${\mathcal F}$, telle que $f_k(x) \to z_{\infty}$, et telle que $A(f_k)=\diag(\lambda_1(k), \ldots , \lambda_n(k))$ v\'erifie $1/\lambda_i(k) \to + \infty$ pour un certain $i \in \{1, \ldots n \}$.  Mais le Fait \ref{fait.dist.borne} impliquerait que $(f_k)$ n'est pas stable en $x$, et une fois de plus, la famille ${\mathcal F}$ ne serait pas \'equicontinue en $x$.

En r\'esum\'e, si l'on appelle $K=K_1K_0$, on a montr\'e qu'il existe deux applications $L_1: {\mathcal F} \to K_0$ et $L_2^{\prime}: {\mathcal F} \to K$ telles que pour tout $f \in {\mathcal F}$, on ait~:
$$ p(f)=L_1(f)A(f)L_2^{\prime}(f).$$

   Appelons ${\mathcal K}.K_0$  le compact de $\hn$ obtenu en prenant les unions des translat\'es \`a droite de ${\mathcal K}$ par des \'el\'ements de $K_0$.  On choisit  ${\mathcal U} $ et ${\mathcal V}$ deux voisinages relativement compacts de $0$ dans $\lien^-$, tels que~:
  \begin{enumerate}
  \item{L'application $\Phi\;: \xi \mapsto \pi_M(\exp(\hx,\xi))$ est d\'efinie et injective sur $\overline{\mathcal U}$, et r\'ealise un diff\'eomorphisme de ${\mathcal U}$ sur un ouvert $U \subset M$ contenant $x$.
}
  \item{Pour tout $\hy \in {\mathcal K}.K$, l'application $\Psi_{\hy}\;:\xi \mapsto \pi_N(\exp(\hy, \xi))$ est injective sur $\overline{\mathcal V}$, et r\'ealise un diff\'eomorphisme de ${\mathcal V}$ sur son image. }
\item{Pour tout $l \in K_2K$,   $(\Ad l).\overline{\mathcal U}$ est inclus dans ${\mathcal V}$.}
  \end{enumerate}
  
%   Enfin,  il existe ${\mathcal U}$ un voisinage de $0$ dans $\lien^-$ tel que pour tout $l \in K_2K$, on ait $(\Ad l).{\mathcal U} \subset {\mathcal V}$, et de plus  
Consid\'erons $(f_k)$ une suite de ${\mathcal F}$.  Nous allons montrer qu'il existe une sous-suite de $(f_k)$ qui converge, au sens de la topologie $C^{\infty}$ sur $U$, vers   une application $f \in {\mathcal C}^{\infty}(U,N)$. Pour all\'eger les expressions, on appellera  $l_k:=A(f_k)L_2^{\prime}(f_k)$, et $\hz_k:=f_k(\hx).L_1(f_k)$.   On peut alors \'ecrire~:
$$ f_k(\hx).l_k^{-1}=\hz_k.$$
Comme $(l_k)$ est \`a valeurs dans le compact $K_2K$, on peut extraire une sous-suite de $(f_k)$ telle que $(\Ad l_k)_{| \lien^-}$ converge vers $L \in \End(\lien^-,\lieg)$, et $\hz_k$ converge vers $\hz \in {\mathcal K}K_0$.  Appelons $f : U \to N$ l'application d\'efinie, pour tout $\xi \in {\mathcal U}$ comme suit~:
$$  f(\pi_M(\exp(\hx,\xi))):=\pi_N(\exp(\hz, L(\xi))).$$
Comme $f=\Psi_{\hz} \circ L \circ \Phi^{-1}$ est une compos\'ee de transformations lisses, on a bien $f \in {\mathcal C}^{\infty}(U,N)$. Nous pouvons \`a pr\'esent achever la d\'emonstration du Th\'eor\`eme \ref{thm.normal} gr\^ace au~:

\begin{lemme}
\label{lem.convergence.lisse}
La convergence de $(f_k)$ vers $f$ est  ${\mathcal C}^{\infty}$  sur l'ouvert $U$.
\end{lemme}
\begin{preuve} on consid\`ere $\hat W$ un voisinage ouvert de $\hz$ dans $\hn$, et ${\mathcal W}$ un voisinage ouvert de $0$ dans $\lien^-$, contenant $L({\mathcal U})$, tels que l'application $\exp: \hat W \times {\mathcal W} \to \hn$ soit d\'efinie.  Il s'agit d'une application lisse.  Pour tout $k \in \N$, on d\'efinit $ \varphi_k : {\mathcal U} \to  \hat W \times {\mathcal W}$ par~:
$$\varphi_k(\xi)=(\hz_k,L_k(\xi)).$$
Comme les $L_k$ sont des applications lin\'eaires, la suite $(\varphi_k)$ converge pour la topologie ${\mathcal C}^{\infty}$ sur ${\mathcal U}$ vers $\varphi_{\infty}: \xi \mapsto (\hz,L(\xi))$.  Maintenant, \'etant donn\'e que sur $U$, on a $f_k=\pi_N \circ \exp \circ \;  \varphi_k \circ \Phi^{-1}$ et $f=\pi_N \circ \exp \circ \varphi_{\infty} \circ \pi_M^{-1}$, on obtient la convergence ${\mathcal C}^{\infty}$ de $(f_k)$ vers $f$. \end{preuve}

%Pour cela, donnons nous $(y_k)$ une suite de l'ouvert $V$ qui converge vers $y \in V$.  On souhaite montrer que $f_k(y_k)$ tend vers $f(y)$.  \'Ecrivons $y_k=\pi_M(\exp(\hx,\xi_k))$ et $y=\pi_M(\exp(\hx,\xi))$, avec $\xi,\xi_k \in {\mathcal V}$.
%On obtient~:
%$$f_k(\pi_M(\exp(\hx, \xi_k)).l_k^{-1})=\pi_N(\exp(\hz_k, (\Ad l_k).\xi_k)).$$
%Le terme de droite tend bien vers $\pi_N(\exp(\hz,L(\xi)))=f(y)$, comme annonc\'e. 
\subsection{Int\'egrabilit\'e de certaines distributions associ\'ees aux  suites stables}
\label{sec.integrabilite}
Nous finissons la Section \ref{sec.dynamique-holonomie} par  un r\'esultat d'int\'egrabilit\'e, qui sera \'egalement un point cl\'e dans la preuve du Th\'eor\`eme \ref{thm.integrabilite-stable}, en Section \ref{sec.stratification}.  Nous consid\'erons ici une suite $(f_k)$ d'immersions conformes de $(M,[g])$ dans $(N,[h])$, et nous supposons que $(f_k)$ tend uniform\'ement sur les compacts de $M$ vers une application $f$.  Le contenu de la Proposition \ref{prop.integrabilite} ci-dessous est que lorsque $f \not \in \Conf(M,N)$, des sous-vari\'et\'es totalement g\'eod\'esiques conformes non triviales apparaissent sur $M$.

\begin{notation}
\label{not.normes}
Jusqu'\`a la fin de l'article, on va munir $\lieg$ d'un produit scalaire euclidien $<\  , \ >_{\lieg}$, et on appellera $|| \ ||_{\lieg}$ la norme associ\'ee.  Si $\xi \in \lieg$ et $r>0$, on appellera  ${\mathcal B(\xi,r)}$, $r>0$ la boule de centre $\xi$ et rayon $r$, relativement \`a la norme $\| \ \|_{\lieg}$.
\end{notation}
Soit $x \in M$. Comme $(f_k)$ tend uniform\'ement vers $f$ sur les compacts de $M$, la suite $(f_k)$ est stable en $x $, et par le Lemme \ref{lem.holonomie-stable}, elle admet en $x$ une suite d'holonomie $(h_k)$ dans $\mbox{A}^+$.  En particulier, l'action adjointe de $\Ad h_k$ sur $\lieg$ est diagonalisable sur $\R$~: il existe des suites positives $\nu_1(k)<\ldots< \nu_m(k)$, ainsi qu'une d\'ecomposition $\lieg=\lieg_1 \oplus \ldots \oplus \lieg_m$, de sorte que $(\Ad h_k)_{| \lieg_j}=\nu_j(k)Id_{\lieg_j}$.  Quitte \`a consid\'erer une sous-suite de $(f_k)$, et \`a remplacer $(h_k)$ par une suite \'equivalente dans $\mbox{P}$, on peut supposer~:
\begin{enumerate}
\item{Il existe $(\hx_k)$ une suite de $\hm$ dans la fibre de $x$, qui converge vers $\hx$, et telle que $\hz_k:=f_k(\hx_k).h_k^{-1}$ converge vers $\hz  \in \hn$.}
\item{Chaque suite $(\nu_j(k))$ admet une limite dans $\R_+ \cup \{ \infty \}$, et $\nu_j(k)=o(\nu_{j+1}(k))$, pour tout $1 \leq j  \leq m-1$.}
\end{enumerate}
Dans ce qui suit, nous allons appeler $l$ le plus grand entier de $\{1, \ldots, m \}$ pour lequel $\sup_{k \in \N}\nu_j(k)< \infty$.  
Pour chaque $j \in \{ 1, \ldots ,m \}$, on appelle $\lieg^{(j)}:=\lieg_1 \oplus \ldots \oplus \lieg_j$.   Soient ${\mathcal U}$ et ${\mathcal V}$ deux voisinages de $0$ dans $\lieg$ tels que $\xi \mapsto \exp(\hx, \xi)$ et $\zeta \mapsto \exp(\hz,\zeta)$ soient des diff\'eomorphismes de ${\mathcal U}$ et ${\mathcal V}$ sur leurs images respectives $\hat U$ et $\hat V$. Quitte \`a restreindre  ${\mathcal V}$, il existe $r_0>0$ tel que ${\mathcal V}  \subset  {\mathcal B}(0,r_0)$ et pour tout $\hz^{\prime} \in \exp(\hz,{\mathcal V})$, l'application $\zeta \mapsto \exp(\hz^{\prime},\zeta)$ soit d\'efinie sur ${\mathcal B}(0,r_0)$, et r\'ealise un diff\'eomorphisme de  ${\mathcal B}(0,r_0)$ sur son image.  On munit $\hn$ d'une m\'etrique riemannienne, et l'on appelle $\hat d$ la distance que cette m\'etrique induit sur l'ouvert $\exp(\hz,{\mathcal V})$.  Localement, sur une vari\'et\'e,  les distances induites par deux m\'etriques riemanniennes diff\'erentes sont \'equivalentes. Ainsi, quitte \`a restreindre ${\mathcal U}$,${\mathcal V}$, et $r_0$, il existe deux constantes strictement positives $C_1$ et $C_2$, telles que si $\hz^{\prime} \in \exp(\hz,{\mathcal V})$ et $\zeta_1,\zeta_2 \in {\mathcal B}(0,r_0)$, on a~:

\begin{equation}
\label{equ.dist.equivalente}
C_1|| \zeta_1-\zeta_2||_{\lieg} \leq \hat d(\exp(\hz^{\prime},\zeta_1),\exp(\hz^{\prime},\zeta_2)) \leq C_2 || \zeta_1-\zeta_2||_{\lieg}.
\end{equation}

%, avec de plus la propri\'et\'e que l'adh\'erence de $L(\lieg^{(l)} \cap {\mathcal B}(0,r))$ soit incluse dans ${\mathcal V}$.

On choisit maintenant ${\mathcal U}$ assez petit pour  qu'il existe $C_3>0$, tel que $\forall j \in \{1, \ldots ,l \}$, $\forall \xi \in {\mathcal U}\cap \lieg^{(j)}$, $\forall k \in \N$, on ait~:
\begin{equation}
\label{equ.inclusion}
|| (\Ad h_k).\xi ||_{\lieg} \leq {C_3\nu_j(k)} \leq \min(\frac{r_0}{2},\frac{C_1}{C_2}\frac{r_0}{4}). 
\end{equation}
Ceci est possible puisque les suites $\nu_j(k)$ ont une limite finie pour tout $j \in \{1,\ldots,l\}$. Pour tout $k \in \N$, on appelle ${\hat S}_j(k)$ la sous-vari\'et\'e $\exp(\hx_k,{\mathcal U}\cap \lieg^{(j)})$.

Notre objectif  est maintenant de montrer la proposition ci-dessous (voir \'egalement le Th\'eor\`eme 2.6 de \cite{zeg.affine} dans le cadre isom\'etrique)~:

\begin{proposition}
\label{prop.integrabilite}
Pour tout $j \in \{1, \ldots ,l \}$, la sous-vari\'et\'e~:
$${\hat S}_j:=\exp(\hx,\lieg^{(j)}\cap {\mathcal B}(0,r_1))$$
 est une feuille int\'egrale de la distribution de $\hm$ d\'efinie par $(\om)^{-1}(\lieg^{(j)})$.  De plus, $S_j:=\pi_M({\hat S}_j)$ est une sous-vari\'et\'e totalement g\'eod\'esique conforme de $M$. 
\end{proposition}

%CITER GHANI.

\begin{preuve} nous voulons montrer que si $\hy=\exp(\hx,\xi)$, o\`u $ \xi \in {\mathcal U}\cap \lieg^{(j)}$, alors $\om(T_{\hy}{\hat S}_j)=\lieg^{(j)}$.  Pour tout $k \in \N$, on pose $\hy_k:=\exp(\hx_k,\xi)$.  On choisit $r_2>0$ assez petit pour que ${\mathcal B}(\xi,r_2) \subset {\mathcal B}(0,r_1)$. Alors ${\hat D}_k:=\exp(\hx_k,{\mathcal B}(\xi,r_2) \cap \lieg^{(j)})$ est inclus dans ${\hat S}_j(k)$ pour tout $k \in \N$.  Pour toute suite $(\hy_k^{\prime})$ telle que $\hy_k^{\prime} \in {\hat D}_k$, on peut \'ecrire $\hy_k^{\prime}=\exp(\hx_k,\xi_k)$, o\`u $\xi_k \in {\mathcal B}(\xi,r_2) \cap \lieg^{(j)}$.  Comme $f_k(\hy_k^{\prime}).h_k^{-1}=\exp(\hz_k,(\Ad h_k). \xi_k)$, on tire des relations (\ref{equ.dist.equivalente}) et (\ref{equ.inclusion})~:
\begin{equation}
\label{equ.ineg1}
 \hat d(f_k(\hy_k^{\prime}).h_k^{-1},\hz_k) \leq C_2C_3 \nu_j(k).
 \end{equation}
 
 On d\'efinit ${\mathcal D}_k \subset {\mathcal U}$ par~: $\exp(\hy_k,{\mathcal D}_k)={\hat D}_k$.  Pour la topologie de Hausdorff, $\hat D_k$ tend vers $\hat D:=\exp(\hx,{\mathcal B}(\xi,r_2) \cap \lieg^{(j)})$, et donc ${\mathcal D}_k$ tend vers ${\mathcal D} \subset \lieg$.  On aura $\om(T_{\hy}{\hat S}_j)=\lieg^{(j)}$ si l'on parvient \`a montrer que ${\mathcal D} \subset \lieg^{(j)}$.  
Dans la suite, nous allons noter $\delta_k:=\sup_{\eta \in {\mathcal D}_k}|| (\Ad h_k).\eta  ||_{\lieg}$.  
 
 \begin{lemme}
 \label{lem.delta.borne}
 Si $j=l$, alors la suite $(\delta_k)$ est born\'ee.
 \end{lemme}
  
 \begin{preuve}
 si ce n'est pas le cas, quitte \`a passer \`a des suites extraites, on va avoir $\delta_k \to \infty$. Par connexit\'e de ${\mathcal D}_k$, il va alors exister $(\eta_k)$ avec $\eta_k \in {\mathcal D}_k$ pour tout $k \in \N$, satisfaisant $||(\Ad h_k).\eta_k||_{\lieg}=r_0$.  Soit $\hy_k^{\prime}:=\exp(\hy_k,(\Ad h_k).\eta_k)$, $\hz_k^{\prime}:=f_k(\hy_k).h_k^{-1}$ et $\hz_k^{\prime \prime}:=f_k(\hy_k^{\prime}).h_k^{-1}$.  Alors l'in\'egalit\'e (\ref{equ.dist.equivalente}) donne $\hat d(\hz_k^{\prime \prime},\hz_k^{\prime}) \geq C_1 r_0$, tandis que (\ref{equ.dist.equivalente}) et  (\ref{equ.inclusion})  donnent~:
 $$\hat d(\hz_k^{\prime \prime},\hz_k) + \hat d(\hz_k^{\prime},\hz_k) \leq 2C_2 \min(\frac{r_0}{2},\frac{C_1}{C_2}\frac{r_0}{4})\leq C_1\frac{r_0}{2}.$$  L'in\'egalit\'e triangulaire conduit alors \`a $C_1 r_0 \leq C_1 \frac{r_0}{2}$~: contradiction. \end{preuve}

Du fait que les limites de $(\nu_{l+1}(k)), \ldots (\nu_m(k))$ sont $+ \infty$, on d\'eduit du Lemme \ref{lem.delta.borne} que pour toute suite $(\eta_k)$ de ${\mathcal D}_k$, les valeurs d'adh\'erences de $(\eta_k)$ doivent \^etre dans $\lieg^{(l)}$.  Donc ${\mathcal D} \subset \lieg^{(l)}$ et ${\hat S}_l$ est bien une feuille int\'egrale de $(\om)^{-1}(\lieg^{(l)})$.  

On va maintenant s'int\'eresser au cas $j<l$.

 \begin{lemme}
 Si $j \in \{1, \ldots ,l-1 \}$, la suite $(\delta_k)$ est en $o(\nu_{j+1}(k))$.
 \end{lemme}

\begin{preuve} si les conclusions du lemme sont inexactes, il existe une constante $C_4>0$ et une suite extraite de $(\delta_k)$, renot\'ee $(\delta_k)$,  satisfaisant $\delta_{k} \geq C_4 \nu_{j+1}(k)$ pour tout $k$. 

Par  connexit\'e de $(\Ad h_k)({\mathcal D}_k)$, il va exister une suite $(\eta_k)$, avec $\eta_k \in {\mathcal D}_{k}$ pour tout $k$, et telle que $|| (\Ad h_k).\eta_k||_{\lieg}=C_4 \nu_{j+1}(k)$.  Quitte \`a choisir $C_4$ assez petit, on aura $C_4 \nu_{j+1}(k) \leq \frac{r_0}{2}$ pour tout $k$.   Si l'on pose $\hz_k^{\prime}:=f_k(\hy_k).h_k^{-1}$ et $\hz_k^{\prime \prime}:=f_k(\exp(\hy_k,\eta_k)).h_k^{-1}=\exp(\hz_k^{\prime},(\Ad h_k).\eta_k)$, on obtient, au vu de (\ref{equ.dist.equivalente})~:
\begin{equation}
\label{equ.maj1}
\hat d(\hz_k^{\prime \prime}, \hz_k^{\prime}) \geq C_1C_4 \nu_{j+1}(k).
\end{equation}
En utilisant (\ref{equ.ineg1}), et en appliquant l'in\'egalit\'e anti-triangulaire \`a $\hz_k,\hz_k^{\prime }$ et $\hz_k^{\prime \prime}$, on aboutit \`a l'in\'egalit\'e  $C_2C_3 \nu_j(k) \geq \hat d(\hz_k^{\prime \prime},\hz_k) \geq C_1C_4\nu_{j+1}(k)-C_2C_3 \nu_j(k)$~:  ceci contredit $\nu_j(k)=o(\nu_{j+1}(k))$. \end{preuve}

Le lemme pr\'ec\'edent nous dit que si $(\eta_k)$ est une suite de $\lieg$ satisfaisant $\eta_k \in {\mathcal D}_k$
   pour tout k, alors $|| (\Ad h_k).\eta_k||_{\lieg}=o(\nu_{j+1}(k))$.  En particulier, on a pour toute suite extraite $\eta_{\sigma(k)}$ que $|| (\Ad h_{\sigma(k)}).\eta_{\sigma(k)}||_{\lieg}=o(\nu_{j+1}(\sigma(k)))$.  Ainsi, toute valeur d'adh\'erence de $(\eta_k)$ doit \^etre dans $\lieg^{(j)}$.  On a bien ${\mathcal D}\ \subset \lieg^{(j)}$, et ${\hat S}_j$ est une feuille int\'egrale de $\lieg^{(j)}$. 
   
   Pour finir, il nous reste \`a v\'erifier que pour $j \in \{1,\ldots,l \}$, $S_j:=\pi_M(\hat S_j)$ est une sous-vari\'et\'e totalement g\'eod\'esique de $M$. On commence par remarquer que $\lieg^{(j)}$ est une sous-alg\`ebre de Lie de $\lieg$. Comme les composantes sur $\lien^-$ et $\liep$ d'un vecteur propre de $\Ad h_k$ sont elles-m\^emes des vecteurs propres de $\Ad h_k$, associ\'es \`a la m\^eme valeur propre, on peut \'ecrire $\lieg^{(j)}=\lien_j^- \oplus \liep_j$ o\`u $\lien_j^- $ et $\liep_j$ sont deux sous-alg\`ebres de Lie de $\lien^-$ et $\liep$ respectivement.  Appelons $\mbox{P}_j$ le sous-groupe connexe de $\mbox{P}$ dont l'alg\`ebre de Lie est $\liep_j$. Soit $\hx \in {\hat S}_j$. La sous vari\'et\'e $\hat \Sigma_j:=\exp(\hx,{\mathcal U}\cap \lien_j^-)$ est transverse aux orbites de l'action \`a droite de $\mbox{P}$, donc \`a celles de $\mbox{P}_j$. Le satur\'e de $\hat \Sigma_j$ par $\mbox{P}_1$ est une sous-vari\'et\'e $\hat M_j \subset {\hat S}_j$ diff\'eomorphe au produit ${\hat \Sigma_j} \times \mbox{P}_j$.  On a $\om(T{\hat M}_j) \subset \lien_j^- \oplus \liep_j$, et il y a en fait \'egalit\'e pour des raisons de dimension. Finalement, $\hat M_j$ est une feuille int\'egrale de $(\om)^{-1}(\lieg^{(j)})$ passant par $\hx$, et $\hat M_j \cap \hat S_j$ est un voisinage ouvert de $\hx$ dans  $\hat S_j$.  Par cons\'equent $\pi_M(\hat M_j \cap \hat S_j)$ est un voisinage ouvert de $x=\pi_M(\hx)$ dans $S_j$, et ce voisinage est, par d\'efinition, un morceau de sous-vari\'et\'e totalement g\'eod\'esique conforme de $M$.  Ceci \'etant valable pour tout $x \in S_j$, on obtient que $S_j$ est elle-m\^eme est  totalement g\'eod\'esique conforme. \end{preuve}

\section{Description de l'adh\'erence de $\Conf(M,N)$ dans $\Cont(M,N)$~: preuve du Th\'eor\`eme \ref{thm.applications.limites}}

\label{sec.description.adherence}
Nous consid\'erons toujours  $(M,[g])$ et $(N,[h])$ deux structures conformes pseudo-riemanniennes de m\^eme signature $(p,q)$, avec $p+q \geq 3$. Nous supposons qu'une  suite d'immersions conformes $f_k : (M,[g]) \to (N,[h])$ converge vers $f \in \Cont(M,N)$, uniform\'ement sur les compacts de $M$.  Par le Th\'eor\`eme \ref{thm.normal}, l'application $f$ est en fait lisse, et la convergence de $(f_k)$ vers $f$ est $C^{\infty}$ sur les compacts. Nous allons d\'ecrire pr\'ecis\'ement l'application $f$, et tirer des cons\'equences g\'eom\'etriques pour $(M,[g])$ lorsque la limite $f$ n'appartient pas \`a $\Conf(M,N)$.

Soit $x \in M$. La convergence uniforme de $(f_k)$ sur les compacts de $M$ assure que $(f_k)$ est stable en $x$. Il existe donc en $x$ une suite d'holonomie $(h_k)$ qui soit dans $\mbox{A}^+ \subset \R_+^* \times \mbox{O}(p,q)$. En particulier,  $h_k=\text{diag}(\lambda_1(k), \ldots , \lambda_n(k))$, et les suites  $\frac{1}{\lambda_i(k)}$ sont born\'ees. On reprend les notations de la Section \ref{sec.integrabilite}, et quitte \`a consid\'erer une suite extraite de $(f_k)$, on fera les m\^emes hypoth\`eses simplificatrices qu'en \ref{sec.integrabilite}~: les transformations $\Ad h_k$ sont simultan\'ement diagonalisables, et l'on peut \'ecrire $\lieg=\lieg_1 \oplus \ldots \oplus \lieg_m$, avec $(\Ad h_k)_{| \lieg_j}=\nu_j(k)Id_{\lieg_j}$, o\`u $\nu_1(k)< \ldots <\nu_m(k)$ sont des suites qui convergent dans $\R_+$, et qui satisfont $\nu_j(k)=o(\nu_{j+1}(k))$. L'entier $l$ est toujours le plus grand entier $j$ de $\{1, \ldots , m\}$ tel que $\lim_{k \to \infty}\nu_j(k)$ soit fini.  Observons que puique la suite $(f_k)$ est stable, $\lien^- \subset \lieg^{(l)}$. La suite  $L_k:=(\Ad h_k)_{|\lien^-}$ est ainsi relativement compacte dans $\text{End}(\lien^-)$.   Quitte \`a consid\'erer une suite extraite de $(f_k)$, qui convergera toujours vers $f$, on peut supposer que $(L_k)$  admet une limite $L \in \End(\lien^-)$, et qu'il existe $\hx_k$ dans la fibre $\pi_M^{-1}(x)$, qui converge vers $\hx$, de sorte que la suite $\hz_k:=f_k(\hx_k).h_k^{-1}$ converge vers $\hz \in \hn$ au-dessus de $f(x)$. Soient  ${\mathcal U}$ et ${\mathcal V}$ comme en section Section \ref{sec.integrabilite}.  Quitte \`a r\'etr\'ecir ces ouverts, on supposera de plus que $\xi \mapsto \pi_M(\exp(\hx,\xi))$ et $\zeta \mapsto \pi_N(\exp(\hz,\zeta))$ sont des diff\'eomorphismes de ${\mathcal U} \cap \lien^-$ et $L({\mathcal U} \cap \lien^-)$ sur leurs images respectives.  On appelle $\hat U=\exp(\hx,{\mathcal U})$ et $U=\pi_M(\exp(\hx,{\mathcal U} \cap \lien^-))$.  Si $y \in U$ s'\'ecrit $y=\pi_M(\exp(\hx,\xi))$, avec $\xi \in {\mathcal U} \cap \lien^-$, on a $f_k(y)=\pi_N(\exp(\hz_k,L_k(\xi)))$, converge vers $\pi_N(\exp(\hz,L(\xi))$.  
On conclut que la restriction de $f$ \`a $U$ est  donn\'ee par l'expression~: 
$$f(\pi_M(\exp(\hx,\xi)))=\pi_N(\exp(\hz,L(\xi))), \text{ pour tout } \xi \in {\mathcal U} \cap \lien^-.$$
%\end{itemize}
%Il est important de noter qu'en passant \`a la limite dans  l'in\'egalit\'e (\ref{equ.inclusion}), on obtient  que $\zeta \mapsto \exp(\hz,\zeta)$ est injective en restriction \`a $L(\mathcal U)$.
%Dans la suite, nous noterons $V:=\pi_M(\exp(\hx,{\mathcal V}))$, et $\hat V:=\exp(\hx,{\mathcal V})$.  Il s'agit d'un voisinage de $x$. Comme $\pi_N$, $\pi_M$, et les applications exponentielles sur $\hm$ et $\hn$ sont lisses, $f$ est lisse sur $V$.  Notre premi\`ere t\^ache va \^etre de prouver~:

Nous alons \`a pr\'esent examiner les diff\'erents cas qui peuvent se pr\'esenter, suivant la nature de l'application $L$.

%\begin{enumerate}
\subsection{ Premier cas~: l'application $L$ appartient \`a $\mbox{GL}(\lien^-)$.}
\label{sec.premier}
Dans ce cas, $L({\mathcal U} \cap \lien^-)$ est un ouvert de $\lien^-$ qui contient l'origine,  et $f$ est un diff\'eomorphisme de $U$ sur son image. Comme la convergence de $(f_k)$ vers $f$ est ${\mathcal C}^{\infty}$, et que les $f_k$ sont conformes, on obtient que $f$ appartient \`a $\Conf(U,N)$.

\subsection{ Second cas~: l'application $L$ est l'application nulle.}
\label{sec.second}
Ce cas interviendra lorsque chaque suite $(\lambda_1(k)), \ldots , (\lambda_n(k))$ tend vers l'infini.  L'application $f$ est alors constante sur $U$.  Ce cas a des cons\'equences g\'eom\'etriques int\'eressantes.  On va en effet montrer la~:

\begin{proposition}
\label{prop.ricci}
Si l'application limite $f$ est constante sur $U$, alors il existe un voisinage $U^{\prime}$ de 
$x$ dans $U$ qui est  conform\'ement Ricci-plat. 
\end{proposition}

\begin{preuve}  il existe un plus grand entier $l_0 \in \{ 1, \ldots , m\}$ tel que pour tout $j \in \{1, \ldots ,l_0 \}$, on ait $\nu_j(k) \to 0$.  Du fait que toutes les suites $(\lambda_1(k)), \ldots , (\lambda_n(k))$ tendent vers l'infini, et que donc leurs inverses tendent vers $0$, on a l'inclusion $\lien^- \subset \lieg^{(l_0)}:=\lieg_1 \oplus \ldots \oplus \lieg_{l_0}$, et \'egalement $\lien^+ \subset \lieg_{l_0+1} \oplus \ldots \oplus \lieg_m$.  En particulier, on conclut que $\lieg^{(l_0)} \subset \lien^- \oplus \lies$  (en effet, on a {\it a priori } $\lieg^{(l_0)} \subset \lien^- \oplus \liez \oplus \lies$ mais la composante selon $\liez$ est triviale puisque comme $h_k \in \mbox{A}^+$, $\Ad h_k$ agit par l'identit\'e sur la composante $\liez$, et ne la contracte donc pas).  On va maintenant pouvoir montrer~:

 \begin{lemme}
 \label{lem.horizontal.integrable}
 Sur $\hat U$, la distribution $(\om)^{-1}(\lien^- \oplus \lies)$ admet une feuille int\'egrale.
 \end{lemme}  

\begin{preuve} d'apr\`es la Proposition \ref{prop.integrabilite}, si $r>0$ est assez petit, alors ${\hat S}_{l_0}:=\exp(\hx,{\mathcal B}(0,r) \cap \lieg^{(l_0)})$ est une sous-vari\'et\'e de $\hat U$, qui est une feuille int\'egrale de $(\om)^{-1}(\lieg^{(l_0)})$.  Quitte \`a choisir $r>0$ encore plus petit, la sous-vari\'et\'e $\hat \Sigma:=\exp(\hx,{\mathcal B}(0,r) \cap \lien^-)$  est transverse aux orbites de l'action \`a droite de $\mbox{S} \subset \mbox{P}$ sur $\hat U$.  Le satur\'e de $\hat \Sigma$ par cette action de $\mbox{S}$ est une sous-vari\'et\'e $\hat M_0 \subset \hat V$, diff\'eomorphe au produit $\hat \Sigma \times \mbox{S}$.  Remarquons que si $\hy_0 \in \hat \Sigma$, alors $\om(T_{\hy_0} \hat \Sigma) \subset \lieg^{(l_0)} \subset \lien^- \oplus \lies$, car $\hat \Sigma \subset {\hat S}_{l_0}$.  Soit maintenant $\hy \in \hm_0$.  On \'ecrit $\hy=\hy_0.p^{-1}$, o\`u $\hy_0 \in \hat \Sigma$ et $p \in \mbox{S}$.  Si ${\mathcal O}$ d\'esigne la $\mbox{S}$-orbite passant par $\hy$, on a le scindement $T_{\hy}\hm_0=T_{\hy}{(\hat \Sigma}.p^{-1})\oplus T_{\hy}{\mathcal O}$.  Les propri\'et\'es d'\'equivariance de la connexion de Cartan conduisent \`a  $\om(T_{\hy}({\hat \Sigma}.p^{-1}))=(\Ad p).\om(T_{\hy_0} {\hat \Sigma})$, d'o\`u $\om(T_{\hy}({\hat \Sigma}.p^{-1})) \subset (\Ad p).\lieg^{(l_0)} \subset\lien^- \oplus \lies$. D'autre part $\om(T_{\hy}{\mathcal O}) \subset \lies$.  Finalement, $\om(T\hat M_0) \subset\lien^- \oplus \lies$, et il y a en fait \'egalit\'e pour des raisons de dimension~:  la sous-vari\'et\'e $\hat M_0$ est une feuille int\'egrale de la distribution $(\om)^{-1}(\lien^- \oplus \lies)$. \end{preuve}

La Proposition \ref{prop.ricci} est maintenant une simple cons\'equence du lemme ci-dessus et de la Proposition \ref{prop.ricci.plate}. \end{preuve}

\subsection{Troisi\`eme cas~: fibration sur une sous-vari\'et\'e totalement d\'eg\'e\-n\'er\'ee}
\label{sec.troisieme}
Si l'on n'est pas dans les deux premiers cas $L$ n'est ni inversible, ni nulle. Par les deux cas trait\'es pr\'ec\'edemment, on v\'erifie que ceci n'est possible que lorsque la signature n'est pas riemannienne, {\it i.e} $p \geq 1$.  Remarquons que puisque $h_k \in \mbox{A}^+ \subset \R_+^* \times \mbox{O}(p,q)$, les suites $(\lambda_i(k))$ s'\'ecrivent $\lambda_i(k)=\sigma_k \alpha_i(k)$, o\`u la matrice $\diag(\alpha_1(k), \ldots , \alpha_n(k))$ est dans $\mbox{O}(p,q)$.  En particulier, les $\alpha_i(k)$ satisfont la relation  $\alpha_i(k)=\frac{1}{\alpha_{n-i+1}(k)}$ pour tout $i \in \{1, \ldots ,p \}$.  Lorsque $q>p$, et  $i \in \{ p+1, \ldots ,n-p \}$, on a $\alpha_i(k)=1$.   Comme $h_k \in \mbox{A}^+$, on a les in\'egalit\'es  $\alpha_1(k)\geq \ldots \geq \alpha_n(k)$, et $\frac{\sigma_k}{\alpha_1(k)} \geq 1$.  On conclut que si $(\sigma_k)$ est born\'ee dans $[1, + \infty[$, il en va de m\^eme pour $(\alpha_1(k))$ et l'on est dans le premier cas ci-dessus, o\`u $L \in \mbox{GL}(\lien^-)$.  Ainsi $(\sigma_k)$ tend vers l'infini.  Comme on n'est pas dans le second cas ci-dessus, o\`u $L$ est constante, il existe $n-p+1 \leq r \leq n$ minimal tel que pour $i \in \{ r ,\ldots ,n \}$, $\sigma_k \alpha_i(k)$ soit  born\'ee dans $[1, +\infty[$.  Soit $(\epsilon_1, \ldots ,\epsilon_n)$ la base de $\lien^-$ d\'efinie par~:
$$ \epsilon_i=  \left( \begin{array}{ccc}
0& -e_i^t.J_{p,q}   &  0\\
  & 0  & e_i \\
  &   & 0
\end{array} \right),
$$
$(e_1, \ldots , e_n)$ \'etant la base canonique de $\R^{p,q}$.  
%Notons que dans la base $(\epsilon_1,\ldots,\epsilon_n)$
Comme l'expression de $(\Ad h_k)_{|\lien^-}$, dans la base $(\epsilon_1, \ldots ,\epsilon_n)$, est $\diag (\frac{\alpha_1(k)}{\sigma_k}, \ldots ,\frac{\alpha_n(k)}{\sigma_k})$, l'image $\mbox{Im } L$ n'est autre que $\text{Vect}(\epsilon_r, \ldots, \epsilon_n)$, et le noyau $\mbox{Ker } L$ est $\text{Vect}(\epsilon_1, \ldots, \epsilon_{r-1})$. 

\subsubsection{Description des fibres de $f_{|U}$}
\label{sec.fibres}
\ \\
De l'expression $f(\pi_M(\exp(\hx,\xi)))=\pi_N(\exp(\hz,L(\xi)))$,  pour tout $ \xi \in {\mathcal U} \cap \lien^-$, et de l'injectivit\'e de $\zeta \mapsto \exp(\hz,\zeta)$ sur $L({\mathcal U} \cap \lien^-)$, on tire que la fibre $F(x):=\{y \in U \ | \ f(y)=f(x) \}$ est donn\'ee par $\pi_M(\exp(\hx, \Ker L \cap {\mathcal U} ))$.   Comme $\exp(\hx, \Ker L \cap {\mathcal U} )$ est transverse aux fibres de $\hm$, $F(x)$ est une sous-vari\'et\'e lisse de $U$.  Le fait qu'elle est totalement g\'eod\'esique conforme r\'esultera du Th\'eor\`eme \ref{thm.integrabilite-stable}.  En effet, dans le cas que nous sommes en train d'\'etudier, l'entier $s$ du  Th\'eor\`eme \ref{thm.integrabilite-stable} est au moins $2$ et la fibre $F(x)$ est une feuille int\'egrale de la distribution ${\mathcal F}_{s-1}$.  La propri\'et\'e qu'ont ces feuilles d'\^etre  totalement g\'eod\'esiques conformes sera montr\'ee en section \ref{preuve2}. 

 Nous avons d\'ej\`a observ\'e que l'espace tangent \`a une sous-vari\'et\'e totalement g\'eod\'esique conforme a une signature constante.  Ici,  la signature des fibres de $f$ est celle de $\Ker L $ relativement \`a  la m\'etrique $\lambda^{p,q}$.  Or $\Ker L $ est l'espace $\text{Vect}(\epsilon_1, \ldots, \epsilon_{r-1})$, o\`u $r>n-p$  a \'et\'e d\'efini ci-dessus.  Dans la base $(\epsilon_1,\ldots,\epsilon_n)$, la forme $\lambda^{p,q}$ (voir la Section \ref{sec.principe-equivalence}) est conforme \`a $2\xi_1\xi_n+ \ldots +2 \xi_p\xi_{q+1}+\Sigma_{p+1}^q \xi_{i}^2$. Le sous-espace $\Ker L $ est donc d\'eg\'en\'er\'e, et son radical est $\text{Vect}(\epsilon_1, \ldots , \epsilon_{n-r+1})$.  Observons que ce radical est de dimension $n-r+1$, qui est exactement la dimension de $\text{Vect}(\epsilon_r, \ldots , \epsilon_n)= \mbox{Im }L$.  Autrement dit, le radical de $F(x)$ a pour dimension le rang de $f$.

%\^eme  fibre que $\exp(\hy,\xi)$, converge vers  $\exp(\hy,\xi)$, et $f_k(\exp(\hy_k,\xi_k)).h_k^{-1}=\exp(f_k(\hy_k).h_k^{-1},L_k(\xi_k))$ tend vers $\exp(\hz,L(\xi))=\hz$.  On en d\'eduit, par la caract\'erisation donn\'ee au Lemme \ref{lem.caracteristique}  que $\exp(\hy,\xi)$ appartient \`a $\hat F(x)$.   
%%Si l'on note $\hz_k^{\prime}:=f_k(\hy).h_k^{-1}$, on a d'une part $\hz_k^{\prime} \to \hz$ par la remarque pr\'ec\'edente, et d'autre part $f_k(\exp(\hy,\xi)).h_k^{-1}=\exp(\hz_k^{\prime},L_k(\xi))$.  On en d\'eduit que $\lim_{k \to \infty}f_k(\exp(\hy,\xi)).h_k^{-1}=\hz$ et donc $\exp(\hy,\xi) \in \hat F(x)$.  
%Comme la dimension de la sous-vari\'et\'e $\hat F(x)$ est celle de $\Ker L$, on conclut que pour ${\mathcal V}^{\prime}$ un petit voisinage de $0$, $\exp(\hy,\Ker L \cap {\mathcal V}^{\prime})$ est un voisinage de $\hy$ dans $\hat F(x)$.  Aussi, $\om_{\hy}(T_{\hy}(\hat F(x)))= \Ker L$.  On voit donc que $F(x)$ est une sous-vari\'et\'e totalement g\'eod\'esique conforme, selon la D\'efinition \ref{def.ss.variete.geodesique}.  

\subsubsection{Description de l'image de $f_{|U}$}
On s'int\'eresse \`a pr\'esent \`a l'image $f(U)$ de l'application $f$.  Il s'agit de $\Sigma=\pi_N(\exp(\hz,L({\mathcal U} \cap \lien^-)))$.  Comme nous l'avons d\'ej\`a signal\'e, $\zeta \mapsto \pi_N(\exp(\hz,\zeta))
$ est un diff\'eomorphisme de $L({\mathcal U} \cap \lien^-)$ sur son image. Ainsi, $\Sigma$ est une sous-vari\'et\'e lisse de $N$. Nous allons en d\'eterminer la signature.

  Soit $y \in U$.  On \'ecrit $y=\pi_M(\exp(\hx,\xi))$, o\`u $\xi \in {\mathcal U} \cap \lien^-$.  Pour $k$ suffisamment grand, on peut \'egalement \'ecrire $y=\pi_M(\exp(\hx_k,\xi_k))$, o\`u $(\xi_k)$ est une suite de $\lien^-$ qui tend vers $\xi$.  On pose $\hy_k:=\exp(\hx_k,\xi_k)$.  Alors $\hz_k^{\prime}:=f_k(\hy_k).h_k^{-1}=\exp(\hx_k,L_k(\xi_k))$ tend vers $\hz^{\prime}:=\exp(\hz,L(\xi))$.  Soit maintenant $0<r_1<<1$ tr\`es petit, tel que $\xi^{\prime} \mapsto \pi_N(\exp(\hy,\xi^{\prime}))$ soit un diff\'eomorphisme de ${\mathcal B}(0,r_1) \cap \lien^-$ sur un voisinage $U_y \subset U$ de $y$.  Si $y^{\prime} \in U_y$, on peut \'ecrire $y^{\prime}=\pi_M(\exp(\hy,\xi^{\prime}))$, avec $\xi^{\prime} \in {\mathcal B}(0,r_1) \cap \lien^-$.  Pour $k$ assez grand, on aura aussi $y^{\prime}=\pi_M(\exp(\hy_k,\xi_k^{\prime}))$, avec $(\xi_k^{\prime})$ une suite de $\lien^-$ qui tend vers $\xi^{\prime}$.  De la relation $f_k(y^{\prime})=\pi_N(\exp(\hz_k^{\prime},L_k(\xi_k^{\prime})))$, on d\'eduit que $f$ est donn\'ee sur $U_y$ par $f(\pi_M(\exp(\hy,\xi^{\prime})))=\pi_N(\exp(\hz^{\prime},L(\xi^{\prime})))$, $\xi^{\prime} \in {\mathcal B}(0,r_1) \cap \lien^-$.  On en d\'eduit que l'espace tangent \`a $f(U)$ en $f(y)$ est $\iota_{\hz^{\prime}}(\mbox{Im }L)$, et ce pour tout $y \in U$.   
La signature de $\Sigma$ est donc celle de $\mbox{Im }L$  relativement \`a $\lambda^{p,q}$.   Comme $\mbox{Im }L$ co\"{\i}ncide avec $\text{Vect}(\epsilon_r, \ldots, \epsilon_n)$,
  donc  est incluse dans $\text{Vect}(\epsilon_{n-p+1}, \ldots , \epsilon_n)$, c' est un sous-espace totalement d\'eg\'en\'er\'e de $\lien^-$, relativement \`a $\lambda^{p,q}$.   La sous-vari\'et\'e $\Sigma$ est \'egalement totalement d\'eg\'en\'er\'ee.  

%totalement d\'eg\'en\'er\'ee pour $\lambda^{p,q}$, 

\subsubsection{Conclusion}  En chaque point $x \in M$, on a d\'ecrit l'allure locale de $f$  selon que l'on \'etait dans le cas d\'ecrit en \ref{sec.premier}, \ref{sec.second} ou \ref{sec.troisieme}.  Ces cas s'excluent mutuellement, et l'ensemble des points conduisant \`a un cas donn\'e est clairement ouvert.  Aussi, par connexit\'e de la vari\'et\'e $M$, l'un des trois cas pr\'ec\'edents d\'ecrit le comportement de $f$ {\it au voisinage de chaque point de $M$}.  
\begin{itemize}
\item{Si l'on est dans le cas \ref{sec.premier}, cela veut dire que $f$ est localement une immersion conforme d'un ouvert de $M$ dans un ouvert de $N$.  Ainsi $f \in \Conf (M,N)$.}
\item {Si l'on est dans le cas \ref{sec.second}, $f$ est localement constante, donc constante par connexit\'e de $M$.  Par ailleurs, la Proposition \ref{prop.ricci} montre chaque point admet un voisinage o\`u une m\'etrique de la classe conforme est Ricci-plate.  La vari\'et\'e $(M,g)$ est localement conform\'ement Ricci-plate dans ce cas.}
\item{Enfin, si l'on est dans le cas \ref{sec.troisieme}, la fibre passant par chaque point est une sous-vari\'et\'e totalement g\'eod\'esique conforme, qui est d\'eg\'en\'er\'ee, et dont le radical a pour dimension le rang de $f$.  Chaque point $x \in M$ admet un voisinage $U$ tel que $f(U)$ soit une sous-vari\'et\'e lisse totalement isotrope de $N$.}
\end{itemize}

%\end{enumerate}

%Tout d'abord, on sait par le Lemme \ref{lem.stabilite-ouverte} qu'il existe une suite stable $h_k \in \mbox{A}^+$ qui soit suite d'holonomie de $(f_k)$ en chaque point de $M$.  On pose comme pr\'ec\'edemmet $L_k:= (\Ad h_k)_{| \lien^-}$, et du fait que $h_k \in \mbox{A}^+$, on a que $(L_k)$ est relativement compacte dans $\End (\lien^-)$.  Quitte \`a consid\'erer une suite extraite de $(f_k)$  (qui convergera toujours vers $f$), on va supposer que $L_k$ converge vers $L \in \End(\lien^-)$. 
\section{Cons\'equences g\'eom\'etriques de la d\'eg\'enerescence en signatures riemanniennes et lorentziennes}
\label{sec.riemann.lorentz}
Consid\'erons  $(M,[g])$ et $(N,[h])$ deux structures pseudo-riemanniennes de m\^eme signature $(p,q)$, $p+q \geq 3$. On aimerait caract\'eriser g\'eom\'etriquement les cas o\`u  $\Conf(M,N)$ n'est pas ferm\'e dans ${\mathcal C}^0(M,N)$, autrement dit les cas o\`u  existe une suite $(f_k)$ d'immersions conformes qui d\'eg\'en\`erent vers une application limite $f:M \to N$ qui n'est pas une immersion conforme.  Nous commen\c{c}ons par \'etablir un r\'esultat technique en Section \ref{sec.uniformite}, puis nous r\'epondrons, au moins localement \`a la question ci-dessus dans le cas des vari\'et\'es riemanniennes et lorentziennes.

\subsection{Uniformit\'e de l'holonomie}
\label{sec.uniformite}
Nous commen\c{c}ons par montrer qu'il existe une m\^eme suite d'holonomie, valable en chaque point de $M$, pour la suite $(f_k)$.

%On commence par montrer que la stabilit\'e (resp. la stabilit\'e forte) est une propri\'et\'e ouverte :

\begin{lemme}
\label{lem.stabilite-ouverte}
On suppose  que $f_k: (M,g) \to (N,h)$ est une suite d'immer\-sions conformes qui converge dans $\Cont(M,N)$ vers une application $f$.  Alors il existe une m\^eme suite stable  $(h_k)$ de $\mbox{P}$,  qui soit une suite d'holonomie de $(f_k)$ en tous les points de $M$.  De plus, quitte \`a consid\'erer une suite extraite de $(f_k)$, il existe dans la fibre $\pi_M^{-1}(x)$ de chaque point $x \in M$, une suite convergente $(\hx_k)$, telle que $f_k(\hx_k).h_k^{-1}$ converge dans $\hn$.
%\end{enumerate}
\end{lemme}

\begin{preuve} soit $x \in M$. Comme $(f_k)$ est suppos\'ee converger vers $f$ uniform\'ement sur les compacts de $M$, $(f_k)$ est stable en $x$. Il existe donc en $x$ une suite d'holonomie $(h_k)$ qui soit dans $\mbox{A}^+ \subset \R_+^* \times \mbox{O}(p,q)$. En particulier,  $h_k=\text{diag}(\lambda_1(k), \ldots , \lambda_n(k))$, $\lambda_1(k) \geq \ldots \lambda_n(k) \geq 1$.  La suite  $L_k:=(\Ad h_k)_{|\lien^-}$ est ainsi relativement compacte dans $\text{End}(\lien^-)$. 
Par d\'efinition d'une suite d'holonomie, il existe une suite $(\hx_k)$ relativement compacte de $\hm$, dans la fibre de $x$,  telle que $f_k(\hx_k).h_k^{-1}$ soit \'egalement relativement compacte dans $\hn$. On se fixe  $U$ un voisinage de $x$ suffisamment petit, de sorte que pour tout $k \in \N$, il existe ${\mathcal U}_k \subset \lien^-$ un voisinage de $0$ pour lequel $\xi \mapsto \pi_M ( \exp(\hx_k,\xi))$ soit un diff\'eomorphisme de ${\mathcal U}_k$ sur $U$. 

Soit $y \in U$. Il existe une suite $\xi_k$ de $\lien^-$, $\xi_k \in {\mathcal U}_k$,  telle que pour tout $k$,    $\pi_M(\exp(\hx_k,\xi_k))=y$. On peut alors \'ecrire :
\begin{equation}
\label{equn0}
 f_k(\exp(\hx_k,\xi_k)).h_k^{-1}=\exp(f_k(\hx_k).h_k^{-1}, (\Ad h_k).\xi_k).
\end{equation}
La suite $\exp(f_k(\hx_k).h_k^{-1}, (\Ad h_k).\xi_k)$ est relativement compacte dans $\hn$, et il s'ensuit que $(h_k)$ est une suite d'holonomie de $(f_k)$ en $y$. 

On d\'efinit une relation d'\'equivalence sur $M$ comme suit~: deux points $x_1$ et $x_2$ de $M$ sont \'equivalents (on note $x_1 \sim x_2$) si et seulement si les suites d'holonomies $(h_k)$ et $(h_k^{\prime})$ de $(f_k)$ en $x_1$ et $x_2$ respectivement sont \'equivalentes dans $\mbox{P}$ (voir section \ref{sec.holonomie-suite}).  L'argument pr\'ec\'edent montre que les classes d'\'equivalence de la relation $\sim$ sont des ouverts de $M$.  Par connexit\'e de $M$, il ne peut donc y avoir qu'une seule classe d'\'equivalence, ce qui prouve le premier point du lemme.

On reprend les notations ci-dessus.  Quitte \`a consid\'erer une suite extraite de $(f_k)$ (et la suite extraite correspondante de $(h_k)$),  on peut supposer que $\hx_k$ tend vers $\hx$ dans la fibre au-dessus de $x$, que $f_k(\hx_k).h_k^{-1}$ converge vers $\hz \in \hn$, et que $L_k$ converge vers $L \in \text{End}(\lien^-)$.  Alors, si $y \in U$, on peut \'ecrire pour tout $k \in \N$, $\pi_M(\exp(\hx_k,\xi_k))=y$ avec $\xi_k \in {\mathcal U}_k$.  Cette fois, la suite $(\xi_k) $ tend vers $\xi \in \lien^-$.  La relation (\ref{equn0}) assure alors que $f_k(\exp(\hx_k,\xi_k)).h_k^{-1}$ converge vers $\exp(\hz,L(\xi))$.  Notons que $\exp(\hx_k,\xi_k)$ est une suite de la fibre $\pi_M^{-1}(y)$ qui converge vers $\hy:=\exp(\hx,\xi)$.
En r\'esum\'e, nous venons de montrer que chaque point $x \in M$ admet un voisinage ouvert $U$ avec la propri\'et\'e~: {\it de toute suite extraite  $(f_{\sigma(k)})$ de $(f_k)$, on peut extraire une nouvelle sous-suite $(f_{\varphi(k)})$, de sorte que pour tout $y \in U$, il existe dans la fibre $\pi_M^{-1}(y)$ une suite convergente $(\hy_k)$, telle que $f_{\varphi(k)}(\hy_k).h_{\varphi(k)}^{-1}$ converge dans $\hn$.}
Comme $M$ est r\'eunion d\'enombrable de compacts, un proc\'ed\'e diagonal standard permet d'obtenir le second point du lemme. \end{preuve}

\subsection{Th\'eor\`eme \ref{thm.ferme.ou.plat}~: le cadre riemannien}
\label{sec.riemann}
Le Th\'eor\`eme \ref{thm.applications.limites} donne les diff\'erentes possibilit\'es pour l'application limite $f$. Dans le cas riemannien, il s'agit n\'ecessairement d'une application constante.  
Quitte \`a consid\'erer une sous-suite de $(f_k)$ (qui convergera toujours vers $f$), nous allons supposer que les conclusions du Lemme \ref{lem.stabilite-ouverte} sont satisfaites.
  Par le Lemme \ref{lem.stabilite-ouverte}, la suite $(f_k)$ admet une m\^eme suite d'holonomie stable $(h_k)$ en tout point de $M$, et comme on est en signature riemannienne, la suite $(h_k)$ est de la forme $\diag(\lambda_k,\ldots,\lambda_k) \in\mbox{A}^+ $, o\`u $1/\lambda_k \to 0$.  Notons que pour tout $\xi \in \lien^-$, on a $(\Ad h_k).\xi=\frac{1}{\lambda_k}.\xi$. Soit $x \in M$, et soit $(\hx_k)$ une suite   de la fibre $\pi_M^{-1}(x)$ qui converge vers $\hx$, telle que $\hz_k:=f_k(\hx_k).h_k^{-1}$ converge vers $\hz \in \hn$.  Comme $\Ad h_k$ agit par l'homoth\'etie de rapport $1/\lambda_k$  (resp. $\lambda_k$) sur $\lien^-$ (resp. sur $\lien^+$), et trivialement sur $\liez \oplus \lies$,  la relation d'\'equivariance sur la courbure (voir (\ref{equ.equivariance.courbure}) en Section \ref{sec.courbure})~:
$$ (\Ad h_k^{-1}).\kappa_{\hz_k}((\Ad h_k).\xi, (\Ad h_k).\eta)=\kappa_{\hx_k}(\xi,\eta)$$
montre que $\kappa_{\hx}(\xi,\eta)=0$ pour tout $\xi,\eta \in \lien^-$.  Ainsi la courbure de la connexion normale de Cartan associ\'ee \`a $(M,[g])$ s'annule~: $(M,g)$ est localement conform\'ement plate.

\begin{remarque}
Dans \cite{ferrand2}, J. Ferrand obtient, en supposant que les applications $(f_k)$ sont injectives, que $(M,g)$ est en fait conform\'ement \'equivalent \`a un ouvert de $\R^n$.  On pourrait retrouver ce r\'esultat avec les m\'ethodes pr\'esent\'ees ci-dessus, en s'inspirant de \cite{frances-ferrand}, Section 6.
\end{remarque}

\subsection{Th\'eor\`eme \ref{thm.ferme.ou.plat}~: le cadre lorentzien}
\label{sec.lorentz}
On suppose d\'esormais que $(M,g)$ et $(N,h)$ sont lorentziennes, de dimension $\geq 3$. 
L\`a encore, par le Lemme \ref{lem.stabilite-ouverte}, la suite $(f_k)$ admet une m\^eme suite d'holonomie stable $(h_k)$ en tout point de $M$, et comme on est en signature lorentzienne, la suite $(h_k)$ est de la forme $h_k=\diag(\sigma_k\lambda_k, \sigma_k, \ldots ,\sigma_k,\frac{\sigma_k}{\lambda_k})$, avec $1 \leq \lambda_k \leq \sigma_k$, et $\sigma_k \to \infty$.  L'action de $\Ad h_k$ sur $\lieg$ est diagonale.  L'espace $\lien^-$ est la somme de trois espaces propres $\lien_1^-$, $\lien_2^-$ et $\lien_3^-$, associ\'es respectivement aux valeurs propres $\frac{1}{\sigma_k\lambda_k}$, $\frac{1}{\sigma_k}$ et $\frac{\lambda_k}{\sigma_k}$.  Sur $\liep=\liez \oplus \lies \oplus \lien^+$, les valeurs propres sont  $\lambda_k,1,\frac{1}{\lambda_k}, \sigma_k \lambda_k, \sigma_k$ et $\frac{\sigma_k}{\lambda_k}$.  On en d\'eduit ais\'ement le~:
\begin{fait}
\label{fait.1}
Si $(\eta_k)$ est une suite de $\liep$ qui converge vers $\eta$, alors il existe une constante $C>0$, qui d\'epend de la suite $(\eta_k)$,  telle que $||(\Ad h_k^{-1}).\eta_k||_{\lieg} \leq C \lambda_k.$
\end{fait}

Soit $x \in M$, et soit $(\hx_k)$ une suite   de la fibre $\pi_M^{-1}(x)$ qui converge vers $\hx$, telle que $\hz_k:=f_k(\hx_k).h_k^{-1}$ converge vers $\hz \in \hn$. 
Du Fait \ref{fait.1}, on d\'eduit le~:
\begin{fait}
\label{fait.2}

\begin{enumerate}
\item{Soit $(\xi,\eta) \in \lien_1^- \times \lien_2^-$;  il existe une constante $C>0$ telle que~:
$$||(\Ad h_k^{-1}).\kappa_{\hz_k}((\Ad h_k).\xi, (\Ad h_k).\eta)||_{\lieg} \leq C \frac{1}{\sigma_k^2}.$$}
\item{Soit $(\xi,\eta) \in \lien_1^- \times \lien_3^-$;  il existe une constante $C>0$ telle que~:
$$||(\Ad h_k^{-1}).\kappa_{\hz_k}((\Ad h_k).\xi, (\Ad h_k).\eta)||_{\lieg}\leq C \frac{\lambda_k}{\sigma_k^2}.$$}
\item{Soit $(\xi,\eta) \in \lien_2^- \times \lien_3^-$;   il existe une constante $C>0$ telle que~:
$$||(\Ad h_k^{-1}).\kappa_{\hz_k}((\Ad h_k).\xi, (\Ad h_k).\eta)||_{\lieg}\leq C \frac{\lambda_k^2}{\sigma_k^2}.$$}
\end{enumerate}
\end{fait}

Le Th\'eor\`eme \ref{thm.applications.limites} donne deux posibilit\'es pour l'application limite $f$~: 
 il s'agit soit d'une application constante, soit d'une submersion lisse sur un segment g\'eod\'esique de lumi\`ere de $(N,h)$.  
 
 \subsubsection{Si $f$ est constante, $(M,g)$ est localement conform\'ement plate}
Dans ce cas, la preuve du Th\'eor\`eme \ref{thm.applications.limites}  montre que les suites $\sigma_k$ et $\lambda_k$ v\'erifient, outre les propri\'et\'es \'enonc\'ees ci-dessus, $\frac{\sigma_k}{\lambda_k} \to 0$.   
R\'e\'ecrivons la relation d'\'equivariance sur la courbure~:
$$ (\Ad h_k^{-1}).\kappa_{\hz_k}((\Ad h_k).\xi, (\Ad h_k).\eta)=\kappa_{\hx_k}(\xi,\eta).$$
Comme les trois suites $\frac{1}{\sigma_k^2}$, $\frac{\lambda_k}{\sigma_k^2}$ et $\frac{\lambda_k^2}{\sigma_k^2}$ tendent  vers $0$, le Fait \ref{fait.2} conduit  \`a $\kappa_{\hx}(\xi,\eta)=0$ pour tout $\xi,\eta \in \lien^-$; la vari\'et\'e $(M,g)$ est localement conform\'ement plate.

\subsubsection{Si $f$ est une submersion sur une g\'eod\'esique de lumi\`ere de $(N,h)$, $(M,g)$ est localement conform\'ement Ricci-plate}
Dans ce cas, la suite $(\frac{\sigma_k}{\lambda_k})$ admet une limite finie dans $]0,1]$.  Nous pouvons alors supposer, quitte \`a multiplier \`a droite $(h_k)$ par une suite convergente de $\mbox{A}^+$, que $\lambda_k=\sigma_k$ pour tout $k \in \N$, et que cette limite est $1$.  L'action adjointe de $\Ad h_k$ sur $\lieg$ est diagonale, avec pour valeurs propres $\nu_1(k)=\frac{1}{\sigma_k^2}, \nu_2(k)=\frac{1}{\sigma_k},\nu_3(k)=1, \nu_4(k)=\sigma_k$ et enfin $\nu_5(k)=\sigma_k^2$.  Les sous-espaces propres associ\'es sont  $\lieg_1, \ldots, \lieg_5$, comme en Section \ref{sec.integrabilite}.  On appelle ${\mathcal G}^+=\lieg_1 \oplus \lieg_2 \oplus \lieg_3$, le sous-espace faiblement stable associ\'e \`a $\Ad h_k$.  La Proposition \ref{prop.integrabilite} assure que pour $r_1>0$ choisi suffisamment petit, ${\hat S}_+:=\exp(\hx, {\mathcal  B}(0,r_1) \cap {\mathcal G}^+)$ est une feuille int\'egrale de $(\om)^{-1}({\mathcal G}^+)$.  Appelons ${\mathcal H}^+$ le sous-espace $\lien^- + \lieg_2$.  Comme ${\mathcal H}^+ \subset {\mathcal G}^+$, la  distribution $\hat{\mathcal H}:=\{(\om)^{-1}({\mathcal H}^+) \ | \ \hy \in {\hat S}_+ \}$ est une distribution de $T{\hat S}_+$.  

\begin{lemme}
\label{lem.integrabilite.2}
La distribution $\hat{\mathcal H} \subset T{\hat S}_+$ est int\'egrable.  
\end{lemme}

\begin{preuve}
Soit $\xi_1, \ldots, \xi_m$ une base de ${\mathcal H}^+$, de sorte que chaque $\xi_i$ est soit dans $\lien^-$, soit dans $\lieg_2\cap \lies$, et $\hat X_1, \ldots , \hat X_m$ les $m$ champs $\om$-constants associ\'es sur $\hm$.  Rappelons que la courbure de la connexion $\om$, \'evalu\'ee sur les $\hat X_i$, est donn\'ee par la formule $d \om(\hat X_i, \hat X_j)+[\om(\hat X_i), \om(\hat X_j)]$, si bien que~:
$$[\xi_i,\xi_j] -\kappa_{\hy}(\xi_i,\xi_j)=\om([\hat X_i, \hat X_j](\hy)), \ \ j=1, \ldots ,m.$$
L'objectif est de montrer que $\hat{\mathcal H}$ est involutive, {\it i.e} pour tout $\hy \in {\hat S}_+$, $\om([\hat X_i, \hat X_j](\hy)) \in {\mathcal H}^+$.  

Si $(\xi_i,\xi_j) \in \lien^- \times (\lieg_2\cap \lies)$, ou $(\xi_i,\xi_j) \in (\lieg_2\cap \lies) \times (\lieg_2\cap \lies)$, alors $\kappa_{\hy}(\xi,\xi)=0$, et donc $\om([\hat X_i, \hat X_j](\hy)) \in {\mathcal H}^+$ puisque $[\lien^-,\lieg_2 \cap \lies] \subset \lien^-$, et $[\lieg_2 \cap \lies,\lieg_2 \cap \lies] \subset \lieg_1 \subset \lien^-$.

Il reste \`a examiner le cas o\`u $\xi_i$ et $\xi_j$ sont tous deux dans $\lien^-$.  On \'ecrit $\hy:=\exp(\hx,\xi)$, o\`u $\xi \in {\mathcal B}(0,r_1) \cap {\mathcal G}^+$.  Du fait que $(\Ad h_k)_{| {\mathcal G}^+}$ converge dans $\End({\mathcal G}^+, \lieg)$, on a que $\hz_k^{\prime}:=f_k(\exp(\hx_k,\xi)).h_k^{-1}=\exp(\hz_k, (\Ad h_k). \xi)$ converge vers $\hz^{\prime} \in \hn$.  Appelons $\hy_k:=\exp(\hx_k,\xi)$.  De la relation d'\'equivariance~:
$$ (\Ad h_k^{-1}).\kappa_{\hz_k^{\prime}}((\Ad h_k).\xi_i, (\Ad h_k).\xi_j)=\kappa_{\hy_k}(\xi_i,\xi_j),$$
et des points $(1)$ et $(2)$ du Fait \ref{fait.2}, on tire que $\kappa_{\hy}(\xi_i,\xi_j) \in \lieg_1$ si $(\xi_i,\xi_j) \in \lien_1^- \times \lien_2^-$, et $\kappa_{\hy}(\xi_i,\xi_j) \in \lieg_1\oplus \lieg_2 \subset {\mathcal H}^+$ si $(\xi_i,\xi_j) \in \lien_1^- \times \lien_3^-$.  On a encore $\om([\hat X_i, \hat X_j](\hy)) \in {\mathcal H}^+$ dans ces deux cas.

Enfin, si $(\xi_i,\xi_j) \in  \lien_1^- \times \lien_3^-$, alors la relation d'\'equivariance pour la courbure s'\'ecrit~:
$$\frac{1}{\sigma_k}\kappa_{\hz_k^{\prime}}(\xi_i,\xi_j)=(\Ad h_k). \kappa_{\hy_k}(\xi_i,\xi_j).$$

Comme $\lieg_2 \cap \lies$ est le seul sous-espace propre contract\'e de $\liep$, on \'etablit facilement le~:
\begin{fait}
\label{fait.3}
Si $(\eta_k)$ est une suite de $\liep$ qui converge vers $\eta$, et si $||(\Ad h_k).\eta_k||_{\lieg} \to 0$, alors $\eta \in \lieg_2$.
\end{fait}

Par le Fait \ref{fait.3}, on conclut que $\kappa_{\hy}(\xi_i,\xi_j) \in \lieg_2 \cap \lies \subset {\mathcal H}^+$, et finalement $\om([\hat X_i, \hat X_j](\hy)) \in {\mathcal H}^+$ lorsque $(\xi_i,\xi_j) \in  \lien_1^- \times \lien_3^-$.  Ceci ach\`eve la preuve du lemme. \end{preuve}

On va en d\'eduire que $(\om)^{-1}(\lien^-\oplus \lies)$ admet une feuille int\'egrale dans $\hm$, passant par $\hx$, et conclure qu'un voisinage de $x$ est conform\'ement Ricci-plat, gr\^ace \`a la Proposition \ref{prop.ricci.plate}.  La preuve est la m\^eme que celle du Lemme \ref{lem.horizontal.integrable}~: on commence par consid\'erer $\hat H$, la feuille int\'egrale de $\hat{\mathcal H}$ passant par $\hx$. Pour $r>0$ assez petit, la sous-vari\'et\'e $\hat \Sigma:=\exp(\hx,{\mathcal B}(0,r) \cap \lien^-)$  est transverse aux orbites de l'action \`a droite de $\mbox{S}$ sur $\hn$.  De plus $\hat \Sigma$ est une sous-vari\'et\'e de ${\hat H}$; en particulier $\om(T \hat \Sigma) \subset {\mathcal H}^+ \subset \lien^- \oplus \lies$.   Le satur\'e de $\hat \Sigma$ par l'action de $\mbox{S}$ est une sous-vari\'et\'e $\hat M_0$, et on a $\om(T \hat M_0) \subset (\Ad \mbox{S}).{\mathcal H}^+ + \lies \subset \lien^- \oplus \lies$.  On a \'egalit\'e pour des raisons de dimension, et finalement $\hm_0$ est une feuille int\'egrale de $(\om)^{-1}(\lien^-\oplus \lies)$.

\section{Le th\'eor\`eme de stratification dynamique \ref{thm.integrabilite-stable}}
\label{sec.stratification}
Dans tout ce qui suit, on va munir la vari\'et\'e pseudo-riemannienne $(N,h)$ d'une m\'etrique riemannienne auxiliaire $\lambda$, qui d\'efinit une distance $d$ sur $N$.  Si $u$ est un vecteur de $TN$, on d\'esignera par $||u||$ la norme de ce vecteur pour la m\'etrique $\lambda$.  Bien entendu, les \'enonc\'es seront ind\'ependants de ce choix d'une m\'etrique auxiliaire.  On consid\`ere une suite d'immersions conformes $f_k : (M,g) \to (N,h)$, et l'on suppose que $(f_k)$ tend vers $f \in \Cont(M,N)$, uniform\'ement sur les compacts de $M$.   La preuve du Th\'eor\`eme \ref{thm.integrabilite-stable} va \^etre l'objet des trois sections \ref{preuve0}, \ref{preuve1}, \ref{preuve2} ci-dessous.

\subsection{D\'efinition des suites $(\mu_j(k))_{k \in \N}$ et des distributions ${\mathcal F}_j$.}
%On suppose que $(f_k)$ est stable en chaque point de $U$. 
\label{preuve0}
On commence par remplacer $(f_k)$ par une suite extraite, renot\'ee $(f_k)$, qui satisfait aux deux conclusions du Lemme \ref{lem.stabilite-ouverte}.  L'une de ces conclusions est l'existence d'une suite $(h_k)$ de $\mbox{A}^+$ qui soit suite d'holonomie de $(f_k)$ en chaque point de $M$.  C'est cette suite d'holonomie qui, apr\`es d'autres extractions \'eventuelles, va d\'eterminer les suites $\mu_1(k), \ldots, \mu_s(k)$ du Th\'eor\`eme \ref{thm.integrabilite-stable}.

\'Ecrivons $h_k=\text{diag}(\lambda_1(k), \ldots , \lambda_n(k))$, avec $\lambda_1(k) \geq \ldots \geq \lambda_n(k) \geq 1$. Il existe $s$ entiers non nuls $n_1, \ldots, n_s$ tels que pour tout $l \in \{0, \ldots, s-1 \}$, et tout couple d'indices $n_{l}+1 \leq i \leq j \leq n_{l+1}$, le quotient $\frac{\lambda_i(k)}{\lambda_j(k)}$ est born\'e dans $[1,+\infty)$  (on a adopt\'e la convention $n_0=0$). Quitte \`a remplacer $(h_k)$ par une suite \'equivalente de $\mbox{P}$, on peut alors supposer que pour tout $l \in \{0, \ldots, s-1 \}$, et tout couple d'indices $n_{l}+1 \leq i \leq j \leq n_{l+1}$, $\lambda_i(k)=\lambda_j(k)$ pour tout $k \in \N$.   Si $j \in \{ 1, \ldots , s \}$, on note  $\mu_j(k)=\frac{1}{\lambda_{n_j}(k)}$. On a $\mu_1(k) \leq \ldots \leq \mu_s(k) \leq 1$, et l'action de $(\Ad h_k)$ sur $\lien^-$ se fait par une transformation diagonale de valeurs propres $\mu_1(k), \ldots , \mu_s(k)$. Quitte \`a extraire encore, on peut supposer  que chaque suite  $(\mu_j(k))$ admet une limite dans $\R^+$, et que $\frac{\mu_{j+1}(k)}{\mu_j(k)} \to \infty$, pour tout $j \in \{1,\ldots,s-1 \}$.  Notons que  les suites $\mu_j(k)$ tendent vers $0$, sauf \'eventuellement pour $j=s$.
%
% si le quotient $\frac{\lambda_i(k)}{\lambda_j(k)}$ est contenu dans un compact de $]0,+\infty[$, alors $\lambda_i(k)=\lambda_j(k)$ pour tout $k \in \N$sont \'equivalentes (en tant que suites de $\R$), elles sont \'egales.  Apr\`es cette simplification, on n'a plus que $s$ suites distinctes, non \'equivalentes,  sur la diagonale de $h_k^{-1}$, o\`u $s \in \{ 1, \ldots , n \}$. On note ces suites $\mu_j(k)$, $j \in \{ 1, \ldots , s \}$, et l'   

On d\'ecompose $\lieg/\liep$ en somme directe :
$$ \lieg/\liep={\mathcal N}_1 \oplus \ldots \oplus {\mathcal N}_s,$$
o\`u $(\Ad h_k)_{|{\mathcal N}_j}=\mu_j(k)Id_{{\mathcal N}_j}$. On pose de plus ${\mathcal N}_0:= \{0 \}$.  Si $n_j$ d\'esigne la dimension de ${\mathcal N}_j$, alors du fait que $\Ad h_k$ agit conform\'ement pour $\lambda^{p,q}$, on a bien que la matrice $$   \left( \begin{array}{ccc}
\mu_1(k)I_{n_1} &  & 0 \\
  & \ddots   &  \\
  0&    & \mu_s(k)I_{n_s}
\end{array} \right)$$ pr\'eserve la classe conforme de $2x_1x_n+\ldots+2x_px_{n-p+1}+\Sigma_{p+1}^q x_i^2.$

Soit  ${\mathcal E}_0= \{ 0 \} \subsetneq {\mathcal E}_{1} \subsetneq \ldots \subsetneq {\mathcal E}_{s-1}$ la filtration d\'efinie par ${\mathcal E}_j:={\mathcal N}_{1} \oplus \ldots \oplus {\mathcal N}_{j-1} \oplus {\mathcal N}_{j}  $ pour $j \in \{1, \ldots, s-1 \}$.
%Il est alors clair qu'un vecteur $\xi \in \lien^-$ non nul est dans ${\mathcal E}_j$, $j=-s, \ldots, s-1$, si et seulement si  :
%$$ ||(\Ad h_k).\xi||_{\lien^-}=\mbox{O}(\mu_{j+1}(k)).$$
%Par ailleurs $\xi$ est dans ${\mathcal E}_j \setminus {\mathcal E}_{j-1}$

%Nous consid\'erons \`a nouveau l'isomorphisme $\iota_{\hy} : \lieg/\liep \to T_yM$ d\'efini en section \ref{sec.principe-equivalence}.  La relation d'\'equivariance (\ref{equ.equivariance}) \'enonc\'ee dans cette section montre que si  $f$ est une immersion conforme de $(M,g)$ dans $(N,h)$, alors :
%\begin{equation}
%\label{equ.image-iota}
%D_xf(\iota_{\hx}(\xi))=\iota_{f(\hx)}(\xi). 
%\end{equation}

Soit  $x \in M$.  Comme $(f_k)$ satisfait aux conclusions du Lemme \ref{lem.stabilite-ouverte}, il existe une suite $(\hx_k)$ dans la fibre $\pi_M^{-1}(x)$ qui tend vers $\hx \in \pi_M^{-1}(x)$, de sorte que $\hz_k:=f_k(\hx_k).h_k^{-1}$ tende vers $\hz \in \hn$ au-dessus de $f(x)$.  On d\'efinit alors~:
$$ {\mathcal F}_j(x):=\iota_{\hx}({\mathcal E}_j), \ j \in \{ 0, \ldots, s-1 \},$$
o\`u $\iota_{\hx} : \lieg/\liep \to T_xM$ d\'esigne l'isomorphisme d\'efini en section \ref{sec.principe-equivalence}.  \`A premi\`ere vue, le sous-espace ${\mathcal F}_j(x)$  semble d\'ependre du point $\hx$, limite de la suite $(\hx_k)$.  En fait, il n'en est rien, comme le montre le~:
\begin{lemme}
\label{lem.independance}
Soient $(\hx_k)$ et $(\hx_k^{\prime})$ deux suites de la fibre $\pi_M^{-1}(x)$, qui convergent vers $\hx$ et $\hx^{\prime}$ respectivement, et telles que $f_k(\hx_k).h_k^{-1}$ et $f_k(\hx_k^{\prime}).h_k^{-1}$ convergent dans $\hn$.  Alors  $\iota_{\hx}({\mathcal E}_j)=\iota_{\hx^{\prime}}({\mathcal E}_j)$, pour tout $ \ j \in \{ 0, \ldots, s-1 \}$.
\end{lemme}

\begin{preuve} les points $\hx$ et $\hx^{\prime}$ \'etant dans la m\^eme fibre, on \'ecrit~: 
$\hx^{\prime}=\hx.p$, avec $p \in \mbox{P}$. D'apr\`es la relation (\ref{equ.equivariance}), le lemme revient \`a  montrer que $(\Ad p).{\mathcal E}_j={\mathcal E}_j$  pour tout $j=1, \ldots ,s-1$. 

On \'ecrit $\hx_k^{\prime}=\hx_k.p_k$, o\`u $(p_k)$ est une suite de $\mbox{P}$ qui tend vers $p$. On a alors :
$$ f_k(\hx_k^{\prime}).h_k^{-1}=f_k(\hx_k).p_kh_k^{-1}=f_k(\hx_k).h_k^{-1}.h_kp_kh_k^{-1}.$$
Comme, par hypoth\`ese, $f_k(\hx_k).h_k^{-1}$ et $f_k(\hx_k^{\prime}).h_k^{-1}$ convergent, la suite  $h_kp_kh_k^{-1}$ doit converger vers $\tilde p \in \mbox{P}$ car l'action de $\mbox{P}$ sur $\hn$ est propre et libre. Supposons par l'absurde qu'il existe $\xi \in {\mathcal N}_{j_0}$ tel que  $(\Ad p).\xi  \not \in {\mathcal E}_{j_0}$.  Pour tout $k$, $(\Ad p_k).\xi=\Sigma_{i>j_0}\xi_{ik}+\xi_k^{\prime}$, avec chaque $\xi_{ik} \in {\mathcal N}_i$, et $\xi_k^{\prime} \in {\mathcal E}_{j_0}$.  Par hypoth\`ese, il existe $m>j_0$ tel que $\xi_{mk} \to \xi_{m\infty} \not = 0$.  Finalement :
$$ (\Ad h_kp_kh_k^{-1}).\xi=\Sigma_{i>j_0} \frac{\mu_i(k)}{\mu_{j_0}(k)}\xi_{ik}+\xi_k^{\prime \prime},$$   o\`u $ \xi_k^{\prime \prime} \in {\mathcal E}_{j_0}.$
Comme $\frac{\mu_m(k)}{\mu_{j_0}(k)} \to \infty$, la composante de $(\Ad h_kp_kh_k^{-1}).\xi$ selon ${\mathcal N}_m$ n'est pas born\'ee, ce qui est absurde puisque $(\Ad h_kp_kh_k^{-1}).\xi \to (\Ad {\tilde l}).\xi$.  \end{preuve}

 \subsection{Caract\'erisation m\'etrique des ${\mathcal F}_j$}
 \label{preuve1}
Pour ce qui suit, on va identifier $\lieg/\liep$ \`a $\lien^-$, via un isomorphisme qui commute \`a l'action adjointe de $\mbox{Z} \times \mbox{S} $.  On verra donc la filtration ${\mathcal E}_{0} \subsetneq \ldots \subsetneq {\mathcal E}_{s-1}$ dans $\lien^-$, et  on consid\'erera $\iota_{\hx}$ comme un isomorphisme entre $\lien^-$ et $T_xM$.  Apr\`es cette identification, la relation~:
\begin{equation}
\label{equn.relation}
 \iota_{\hx.p^{-1}}((\Ad p).{\xi})=\iota_{\hx}(\xi)
\end{equation}  
est encore valable  lorsque $\xi \in \lien^-$, et $p \in \mbox{Z} \times \mbox{S}  \subset \mbox{P}$.  Si $z \in N$, et $\hz \in \hn$ est dans la fibre au-dessus de $z$, on notera \'egalement $\iota_{\hz}$ l'identification entre $\lien^-$ et $T_zN$.  Rappelons que si  $f$ est une immersion conforme de $(M,g)$ dans $(N,h)$, alors~:
\begin{equation}
\label{equ.image-iota}
D_xf(\iota_{\hx}(\xi))=\iota_{f(\hx)}(\xi). 
\end{equation}

On reprend les notations de la section pr\'ec\'edente~: $(\hx_k)$ est une suite de $\hm$ qui tend vers $\hx$ de sorte que $\hz_k:=f_k(\hx_k).h_k^{-1}$ tende vers $\hz \in \hn$ au-dessus de $z:=f(x)$.

Par ailleurs, comme la suite $(\hz_k)$ est contenue dans un compact de $\hn$, il existe deux constantes  strictement positives $c_1$ et $c_2$  telles que pour tout $k \in \N$, et $\xi \in \lien^-$, on ait~:
\begin{equation}
\label{bi-lipschitz}
c_1|| \xi ||_{\lieg} \leq ||\iota_{\hz_k}(\xi) || \leq c_2 || \xi ||_{\lieg}.
\end{equation}

On consid\`ere pour commencer $u \in T_xM$ non nul, et une suite $(u_k)$ de $T_xM$ qui tend vers $u$. On \'ecrit pour tout $k \in \N$, $u_k=\iota_{\hx_k}(\xi_k)$, avec $\xi_k \in \lien^-$.  La suite $\xi_k$ tend vers $\xi:=\iota_{\hx}^{-1}(u)$. Par d\'efinition des distributions ${\mathcal F}_j(x)$, le  vecteur $u$ appartient  \`a $ T_xM \setminus {\mathcal F}_{s-1}(x)$ si et seulement si $\xi \in \lien^- \setminus {\mathcal E}_{s-1}$.  En d\'ecomposant chaque $\xi_k$ selon la somme ${\mathcal N}_1 \oplus \ldots \oplus {\mathcal N}_s$, on voit ais\'ement que ceci est \'equivalent \`a~:
\begin{equation}
\label{equ.1}
||(\Ad h_k).\xi_k ||_{\lieg} \sim c \mu_s(k),
\end{equation}
 pour un certain r\'eel $c>0$ (qui d\'epend de la suite $(\xi_k)$). Des relations (\ref{equn.relation}) et (\ref{equ.image-iota}), on tire par ailleurs~:
\begin{equation}
\label{equ.image-bis}
D_xf_k(u_k)=\iota_{\hz_k}((\Ad h_k).\xi_k).
\end{equation}
Ainsi, la relation (\ref{bi-lipschitz}), jointe \`a (\ref{equ.1}) et (\ref{equ.image-bis}), montre que $u \in T_xM \setminus {\mathcal F}_{s-1}(x)$ si et seulement si $ ||D_xf_k(u_k)||=\Theta(\mu_s(k))  $.  Ceci prouve le point $1,(a)$ du Th\'eor\`eme \ref{thm.integrabilite-stable}.

%Avec les m\^emes notations que ci-dessus, le vecteur $u$ appartient \`a $T_xM \setminus {\mathcal F}_{j-1}$, $j \in \{2, \ldots , s-1 \}$, si et seulement si  $\xi \in \lien^- \setminus {\mathcal E}_{j-1}$.  Ceci \'equivaut au fait que l'une des composantes de $\xi$ selon ${\mathcal N}_j, \ldots ,{\mathcal N}_s$ est non nulle, soit encore $||(\Ad h_k).\xi_k||_{\lieg} \geq c \mu_j(k)$ lorsque $k$ est suffisamment grand, pour une constante $c>0$.  L\`a encore, les relations   (\ref{bi-lipschitz}),  (\ref{equ.1}) et (\ref{equ.image-bis}), montrent que $u \in T_xM \setminus {\mathcal F}_{j-1}(x)$ si et seulement si il existe une  constante $C>0$ telle que pour $k$ assez grand, $ ||D_xf_k(u_k)|| \geq C \mu_j(k)$.  Ceci prouve le point $1,(b)$.

%

Soit \`a pr\'esent $u \in {\mathcal F}_j(x) \setminus {\mathcal F}_{j-1}(x) $, $j \in \{1,\ldots,s-1 \}$.  On peut \'ecrire $u=\iota_{\hx}(\xi)$ o\`u $\xi \in {\mathcal E}_{j}\setminus {\mathcal E}_{j-1}$. Soit $(u_k)$ une suite de $T_xM$ qui converge vers $u$, et on \'ecrit  pour tout $k \in \N$~: $u_k:=\iota_{\hx_k}(\xi_k)$, o\`u $\xi_k$ est dans $\lien^-$ et converge vers $\xi$. On d\'ecompose $\xi=\xi^{(1)}+ \ldots + \xi^{(s)}$ selon la somme directe ${\mathcal N}_1 \oplus \ldots \oplus {\mathcal N}_s$.  On \'ecrit \'egalement $\xi_k=\xi_k^{(1)}+ \ldots + \xi_k^{(s)}$, o\`u  $\xi_k^{(j)}$ appartient \`a ${\mathcal N}_j$. La suite $(\xi_k^{(j)})$ a une limite non nulle, et par cons\'equent, on tire de l'expression $(\Ad h_k).\xi_k=\mu_1(k)\xi_k^{(1)}+ \ldots + \mu_s(k)\xi_k^{(s)}$ que~:

\begin{equation}
\label{equ.6}
 \mu_j(k) =O( || (\Ad h_k).\xi_k ||_{\lieg}).
\end{equation}

Les relations (\ref{equ.image-bis}) et (\ref{bi-lipschitz})  conduisent \`a~:
$$\mu_j(k) =O( ||D_xf_k(u_k)||),$$
 et le point $(1),(b),(i)$ est satisfait. 
 
 Consid\'erons \`a pr\'esent  la suite $u_k:=\iota_{\hx_k}(\xi)$.  Il s'agit d'une suite de $T_xM$ qui converge vers $u$.  On a l'\'equivalence $|| (\Ad h_k).\xi ||_{\lieg} \sim c\mu_j(k)$, o\`u $c>0$.  Les relations (\ref{bi-lipschitz}) et (\ref{equ.image-bis}) conduisent donc \`a~:
$$ ||D_xf_k(u_k)||=\Theta(\mu_j(k) ),$$ 
et le  point $(1),(b),(ii)$ est satisfait.
 
R\'eciproquement, consid\'erons un vecteur non nul $u \in T_xM$, satisfaisant aux conditions  $(1),(b),(i)$ et $(1),(b),(ii)$ pour un certain indice $j_0 \in \{0,\ldots,s-1 \}$.  
Par le point $(1),(a)$, $u$ n'appartient pas \`a $T_xM \setminus {\mathcal F}_{s-1}(x)$.  Par cons\'equent, $u$ appartient \`a un certain ${\mathcal F}_{j_1}(x) \setminus {\mathcal F}_{j_1-1}(x)$, $j_1 \in \{ 0,\ldots,s-1 \}$.  Or, \`a cause de la  propri\'et\'e $\mu_l(k)=o(\mu_{l+1}(k)), \ \forall l \in \{0,\ldots,s\}$, les points $(1),(b),(i)$ et $(1),(b),(ii)$ ne peuvent pas \^etre  v\'erifi\'es pour deux indices distincts $j_0 \not = j_1$.  On conclut que $j_0=j_1$.   On obtient ainsi l'\'equivalence   $(1),(b)$ du Th\'eor\`eme \ref{thm.integrabilite-stable}.

%De m\^eme, soit $(u_k)$ une suite de $T_xM$ qui converge vers $u$, et $\xi_k:=\iota_{\hx_k}^{-1}(u_k)$.  La suite $(\xi_k)$ converge vers $\xi=\iota_{\hx}^{-1}(u)$.  L'\'equivalence $|| D_xf_k(u_k)|| \sim C_u\mu_s(k)$ a lieu pour toute suite $(u_k)$ comme ci-dessus si et seulement si l'\'equivalence $||(\Ad h_k).\xi_k||_{\lien^-} \sim C_u \mu_s(k)$   a lieu pour toute suite $(\xi_k)$ convergeant vers $\xi$.  Ceci caract\'erise  $\xi \in {\mathcal E}_s \setminus {\mathcal E}_{s-1}$, {\it i.e}   $u \in TM \setminus {\mathcal F}_{s-1}$.   
%D'autre part, lorsque $j \in \{1, \ldots , s-1 \}$,   $|| D_xf_k(u_k)|| \sim C \mu_j(k)$ a lieu pour une suite $(u_k)$ comme ci-dessus si et seulement si l'\'equivalence $||(\Ad h_k).\xi_k||_{\lien^-} \sim C \mu_j(k)$   a lieu pour une suite $(\xi_k)$ convergeant vers $\xi$, et ceci caract\'erise $\xi \in {\mathcal E}_j \setminus {\mathcal E}_{j-1}$, {\it i.e } $u \in {\mathcal F}_{j} \setminus {\mathcal F}_{j-1}$.   Ceci ach\`eve la preuve du point $(2)$ du th\'eor\`eme.

\subsection{Int\'egrabilit\'e des distributions ${\mathcal F}_j$ et caract\'erisation m\'etrique des feuilles locales}
\label{preuve2}
Il nous reste \`a prouver l'int\'egrabilit\'e des ${\mathcal F}_j$, le fait que les feuilles $F_j$ sont totalement g\'eod\'esiques conformes, ainsi que la caract\'erisation m\'etrique locale.  Nous gardons les m\^emes notations que ci-dessus : pour $x \in M$, on choisit une suite $(\hx_k)$ de $\hm$ dans la fibre au-dessus de $x$, qui converge vers $\hx \in \hm$, et telle que $\hz_k:=f_k(\hx_k).h_k^{-1}$ converge vers $\hz$.

 On d\'efinit $U_x:=\pi_M(\exp(\hx,{\mathcal B}(0,r_1)))$, o\`u $r_1>0$.
  Quitte \`a choisir $r_1$ assez petit, $U_x$ est un ouvert contenant $x$, et $\xi \mapsto \pi_M(\exp(\hx,\xi))$ est un diff\'eomorphisme de ${\mathcal B}(0,r_1) \cap \lien^-$ sur $U_x$. Pour $j=1, \ldots, s-1$, on d\'efinit $N_j(x):=\pi_M(\exp(\hx,{\mathcal E}_j \cap {\mathcal B}(0,r_1)))$.  
  
 \begin{lemme}
\label{lem.integrable}
 Pour $j \in \{1, \ldots,s-1  \}$,  $N_j(x)$ est une sous-vari\'et\'e lisse, totalement g\'eod\'esique conforme, contenue dans $U_x$, et est partout tangente \`a ${\mathcal F}_j$. 
 \end{lemme}   
 
 \begin{preuve} comme en section \ref{sec.integrabilite},  en consid\'erant \'eventuellement une suite extraite de $(f_k)$, on a une d\'ecomposition $\lieg=\lieg_1 \oplus \ldots \oplus \lieg_m$, de sorte que $(\Ad h_k)_{| \lieg_j}=\nu_j(k)Id_{\lieg_j}$, o\`u  chaque suite $(\nu_j(k))$ admet une limite dans $\R_+ \cup \{ \infty \}$, et $\nu_j(k)=o(\nu_{j+1}(k))$, pour tout $1 \leq j  \leq m-1$.  Il existe des indices $i_1 \leq \ldots \leq i_s$ tels que $\mu_j(k)=\nu_{i_j}(k)$ pour $j \in \{ 1, \ldots ,s \}$.
Si l'on se fixe un indice $j$, alors la Proposition \ref{prop.integrabilite} affirme que quitte \`a restreindre $r_1$, la sous-vari\'et\'e $\hat S_j=\exp(\hx,\lieg^{(i_j)}\cap {\mathcal B}(0,r_1))$ est une feuille int\'egrale de la distribution $(\om)^{-1}(\lieg^{(i_j)})$.  Maintenant, $\exp(\hx,{\mathcal E}_j \cap {\mathcal B}(0,r_1)) \subset \hat S_j$ est transverse aux fibres de $\pi_M$, donc $N_j(x)$ est une sous-vari\'et\'e lisse de $M$. Par ailleurs, on montre comme dans la  preuve de la Proposition \ref{prop.integrabilite}, que $\lieg^{(i_j)}$ est une sous-alg\`ebre de $\lieg$, qui s'\'ecrit $\lieg^{(i_j)}=\lien_{i_j}^- \oplus \liep_{i_j}$, o\`u $\lien_{i_j}^- $ et $ \liep_{i_j}$ sont deux sous-alg\`ebres de $\lien^-$ et $\liep$ respectivement. Par ailleurs, toujours comme dans la  preuve de la Proposition \ref{prop.integrabilite}, on montre que $N_j(x)=\pi_M(\hat S_j)$. Il s'ensuit que $N_j(x)$ est totalement g\'eod\'esique conforme.

 Nous allons maintenant montrer que $N_j(x)$ est tangente \`a la distribution ${\mathcal F}_j$. Soit $y \in N_j(x)$.  On \'ecrit $y=\pi_M(\hy)$ o\`u $\hy:=\exp(\hx,\xi)$, $\xi \in {\mathcal E}_j \cap {\mathcal B}(0,r_1)$.  Alors  $T_yN_j(x)=D_{\hy}\pi_M(T_{\hy}{\hat S}_j)=\iota_{\hy}(\lieg^{(i_j)} \cap \lien^-)=\iota_{\hy}({\mathcal E}_j)$.  Rappelons maintenant que l'on a extrait une sous-suite de $(f_k)$, de sorte que $(\Ad h_k)_{| \lien^-}$ converge vers $L \in \End (\lien^-)$.  Il s'ensuit que $f_k(\exp(\hx_k,\xi)).h_k^{-1}$ converge vers $\exp(\hz,L(\xi))$. Comme $\exp(\hx_k,\xi)$ tend vers $\hy$, on d\'eduit du Lemme \ref{lem.independance} que $\iota_{\hy}({\mathcal E}_j)={\mathcal F}_j(y)$. \end{preuve}
 
 On vient de montrer que la distribution ${\mathcal F}_j$ d\'efinit un feuilletage pour $j = \{1, \ldots ,s-1 \}$, et que $N_j(x)$ est la feuille passant par $x$ dans $U_x$~: on la note d\'esormais $F_j^{loc}(x)$. 
 
\subsubsection{Caract\'erisation m\'etrique des feuilles locales} 
 
Les distances induites par deux m\'etriques riemanniennes sur une vari\'et\'e sont localement \'equiva\-lentes; on a par cons\'equent~:
\begin{lemme}
\label{lem.distance-equivalente}
Il existe un voisinage compact $K$ de $\hz$ dans $\hn$, et des r\'eels strictement positifs $c_1,c_2$ et $r_0$, tels que si $\hz^{\prime} \in K$ et $\xi_1,\xi_2 \in {\mathcal B}(0,r_0) \cap \lien^-$ :
$$ c_1 ||\xi_1-\xi_2 ||_{\lieg} \leq d(\pi_N(\exp(\hz^{\prime},\xi_1)),\pi_N(\exp(\hz^{\prime},\xi_2))) \leq c_2 ||\xi_1-\xi_2 ||_{\lieg}$$
\end{lemme}

On supposera par la suite que l'on a choisi $r_0,r_1>0$ assez petits pour que $\zeta \mapsto \exp(\hz^{\prime},\zeta)$ soit un diff\'eomorphisme de ${\mathcal B}(0,r_0) \cap \lien^-$ sur son image pour tout $\hz^{\prime} \in K$, et que l'on ait de plus, $\forall k \in {\bf N}$~:
$$(\Ad h_k).({\mathcal B}(0,r_1) \cap \lien^-) \subset {\mathcal B}(0,r_0) \cap \lien^-.$$  

On commence par montrer le point $2,(a)$ du Th\'eor\`eme \ref{thm.integrabilite-stable}. Soit $y \in U_x$, et $(y_k)$ une suite de $U_x$ qui tend vers $y$. Alors $y=\exp(\hx,\xi)$, avec $\xi \in {\mathcal B}(0,r_1) \cap \lien^-$, et  pour tout $k$, on peut \'ecrire $y_k=\exp(\hx_k,\xi_k)$, avec $\xi_k \to \xi$. Le point $y$  appartient \`a $U_x \setminus F_{s-1}^{loc}(x)$ si et seulement si $\xi \in \lien^- \setminus {\mathcal E}_{s-1}$.  Par (\ref{equ.1}) et le Lemme \ref{lem.distance-equivalente}, ceci \'equivaut \`a~:
$$ d(f_k(x),f_k(y_k))=\Theta( \mu_s(k)).$$

Soit maintenant $y \in F_j^{loc}(x) \setminus F_{j-1}^{loc}(x)$, $j \in \{1,\ldots,s-1 \}$.  
Alors $y=\pi_M(\exp(\hx,\xi))$, o\`u $\xi \in \{ {\mathcal E}_{j} \setminus {\mathcal E}_{j-1}\} \cap {\mathcal B}(0,r_1)$.  Soit $(y_k)$ une suite de $U_x$ qui converge vers $y$, et on \'ecrit comme ci-dessus  $y_k=\exp(\hx_k,\xi_k)$, avec $\xi_k \to \xi$.  Nous avons d\'ej\`a vu que $\mu_j(k)=O(|| (\Ad h_k).\xi_k  ||_{\lieg})$, ce qui, au vu du Lemme \ref{lem.distance-equivalente}, conduit \`a~:
$$ \mu_j(k)=O(d(f_k(x),f_k(y_k))).$$
Le point $(2),(b),(i)$ est donc satisfait.

Consid\'erons la suite de points $y_k:=\pi_M(\exp(\hx_k,\xi))$, qui  tend vers $y$.
De la relation  $|| (\Ad h_k).\xi ||_{\lieg} \sim c\mu_j(k)$, o\`u $c>0$, et du Lemme \ref{lem.distance-equivalente}, on d\'eduit~: 
\begin{equation}
\label{equ.7}
 d(f_k(x),f_k(y_k)) =\Theta( \mu_j(k)),
\end{equation}
et le point $(2),(b),(ii)$ est satisfait.
% Les relations (\ref{bi-lipschitz}) et (\ref{equ.image-bis}) conduisent donc \`a $C_1 \mu_j(k) \leq ||D_xf_k(u_k)|| \leq C_2\mu_j(k)$, avec $C1,C_2>0$.

R\'eciproquement, soit  $y \in U_x \setminus \{x \}$ satisfaisant  les points $(2),(b),(i)$ et $(2),(b),(ii)$ pour un certain indice $j_0 \in \{0,\ldots s-1 \}$.  Par le $(2),(a)$, $y \not \in U_x \setminus F_{s-1}^{loc}(x)$.  Donc il existe $j_1 \in \{0,\ldots,s-1 \}$ tel que $y \in {F}_{j_1}^{loc}(x) \setminus {F}_{j_1-1}^{loc}(x)$.  Les points $(2),(b),(i)$ et $(2),(b),(ii)$  sont donc \'egalement v\'erifi\'es pour l'indice $j_1$, ce qui implique $j_1=j_0$.     Ceci ach\`eve donc la preuve du  point   $(2),(b)$ du Th\'eor\`eme \ref{thm.integrabilite-stable}.

\section{Appendice~: unicit\'e des feuilles locales   $F_j^{loc}(x)$ et des suites $\mu_j(k)$ dans le Th\'eor\`eme \ref{thm.integrabilite-stable}}
\label{sec.unicite} 
Soit $x$ un point de $M$. Le Th\'eor\`eme \ref{thm.integrabilite-stable} assure l'existence d'une stratification dynamique sur un voisinage $U_x$ de $x$, d\'efinie par des sous-vari\'et\'es ${ F}_0^{loc}(x)=\{x \} \subsetneq {F}_1^{loc}(x) \subsetneq \ldots \subsetneq F_{s-1}^{loc}(x) \subsetneq U_x$, et des suites $\mu_1(k), \ldots \mu_s(k)$, qui satisfont aux conclusions $(2),(a),(b)$ du th\'eor\`eme.  Nous allons prouver l'unicit\'e de cette stratification, au sens o\`u si ${H}_0^{loc}:=\{x \} \subsetneq H_1^{loc}(x) \subsetneq \ldots \subsetneq H_{r-1}^{loc}(x) \subsetneq U_x$ est une autre famille de sous-vari\'et\'es, et si $\eta_1(k), \ldots , \eta_r(k)$ sont d'autres suites, qui v\'erifient  les conclusions $2,(a),(b)$ du Th\'eor\`eme  \ref{thm.integrabilite-stable}, alors $r=s$, $H_j^{loc}(x)=F_j^{loc}(x)$, et $\eta_j(k)=\Theta(\mu_j(k))$, pour tout $j=1, \ldots , s$.  
%
%Si l'un des entiers $r$ ou $s$ vaut $1$, la propri\'et\'e est claire.  Nous supposerons donc dans la suite que $r$ et $s$ sont sup\'erieurs ou  \'egaux \`a $2$. 

Commen\c{c}ons par choisir $y$ et $y^{\prime}$ non nuls, diff\'erents de $x$, dans ${F}_1^{loc}(x)$ et $H_1^{loc}(x)$ respectivement.  Le point $(2), (b),(ii)$ assure qu'il existe une suite $(y_k)$ qui tend vers $y$, telle que $d(f_k(x),f_k(y_k))=\Theta(\mu_1(k))$, et il existe une suite $(y_k^{\prime})$ qui tend vers $y^{\prime}$ telle que $d(f_k(x),f_k(y_k^{\prime})) = \Theta( \eta_1(k))$.    Par le point $(2),(b),(i)$, on doit avoir $\eta_1(k)=\mbox{O}(d(f_k(x),f_k(y_k)))$ et $\mu_1(k)=\mbox{O}(d(f_k(x),f_k(y_k^{\prime})))$, d'o\`u il ressort que $\eta_1(k)=\mbox{O}(\mu_1(k))$ et  $\mu_1(k)=\mbox{O}(\eta_1(k))$.  On aboutit \`a $\eta_1(k)=\Theta(\mu_1(k))$.  Cette information, jointe au $2, (b)$ assure que $y \in H_1^{loc}(x)$ et $y^{\prime} \in {F}_1^{loc}(x)$.  Ainsi $H_1^{loc}(x)={F}_1^{loc}(x)$.

Supposons que pour $j \in \{2,\ldots, \min(r,s)-1 \}$, on ait prouv\'e ${F}_{j-1}^{loc}(x)={ H}_{j-1}^{loc}(x)$.  Choisissons $y \in {F}_{j}^{loc}(x) \setminus {F}_{j-1}^{loc}(x)$, et $y^{\prime} \in {H}_{j}^{loc}(x) \setminus {H}_{j-1}^{loc}(x)$.  Par $(2), (b), (ii)$, il existe $(y_k)$ qui converge vers $y$ et telle que $d(f_k(x),f_k(y_k))=\Theta(\mu_j(k))$.  Par ailleurs, $y \in (U_x \setminus  {H}_{r-1}^{loc}(x)) \bigcup_{i=j}^{r-1} {H}_{i}^{loc}(x)\setminus {H}_{i-1}^{loc}(x)$.  Les points $(2),(a)$ et $(2),(b),(i)$ assurent que $\eta_j(k)=\mbox{O}(d(f_k(x),f_k(y_k)))$.  On conclut que $\eta_j(k)=\mbox{O}(\mu_j(k))$.  Le m\^eme raisonnement appliqu\'e \`a $y^{\prime}$  donne que $\mu_j(k)=\mbox{O}(\eta_j(k))$, et finalement $\eta_j(k)=\Theta(\mu_j(k))$.  Par le point $(2),(b)$, on obtient  $x \in {H}_{j}^{loc}(x) \setminus { H}_{j-1}^{loc}(x)$ et  $y^{\prime} \in {F}_{j}^{loc}(x)\setminus {F}_{j-1}^{loc}(x)$. On conclut finalement que ${F}_{j}^{loc}(x)={ H}_{j}^{loc}(x)$.  

Supposons que $s=\min(r,s)$.  Alors, le raisonnement pr\'ec\'edent permet d'obtenir par induction que pour tout $j \leq s-1$, $\mu_j(k)=\Theta(\eta_j(k))$ et $F_j^{loc}(x)=H_j^{loc}(x)$.  Le point $2,(a)$ du th\'eor\`eme affirme alors que pour tout point $y \in U_x \setminus F_{s-1}^{loc}(x)$ (et donc pour tout $y \in U_x \setminus {H}_{s-1}^{loc}(x)$), et toute suite $(y_k)$ qui converge vers $y$, on a $d(f_k(x),f_k(y_k)) = \Theta(\mu_s(k))$.  On doit ainsi avoir $\eta_j(k)=\Theta(\mu_s(k))$ pour tout $j=s, \ldots, r$.  Comme $\eta_j(k)=o(\eta_{j+1}(k))$, on obtient que $r=s$, et $\eta_s(k)=\Theta(\mu_s(k))$.  Ceci ach\`eve la preuve de l'unicit\'e.

{\bf Remerciements : } je souhaiterais remercier chaleureusement Karin Melnick pour d'int\'eressantes conversations sur le sujet. Par ailleurs, ce travail a b\'en\'efici\'e du  soutien de l'ANR {\sc Geodycos}.

\end{document}